\documentclass{article}
\usepackage{times}
\usepackage{amsfonts}
\usepackage{amssymb}
\usepackage{amsmath}
\usepackage{amscd}
\usepackage{url}
\usepackage{graphicx}
\usepackage[shortlabels]{enumitem}

\def\mcolon{\! : \!}

\newcommand{\qed}{\hfill\mbox{$\Box$}}

\newfont{\msbm}{msbm10 scaled\magstephalf}
\usepackage[linecolor=black,backgroundcolor=mygray,colorinlistoftodos,prependcaption,textsize=small]{todonotes}
\usepackage{xargs}
\usepackage{xcolor}
\definecolor{mygray}{gray}{0.85}

\newcommand\myrestriction{\mathord\restriction}

\def\rg{\rm rg}

\newtheorem{theorem}{Theorem}[section]

\newtheorem{lemma}[theorem]{Lemma}

\newtheorem{definition}[theorem]{Definition}

\newtheorem{claim}[theorem]{Claim}
\newtheorem{fact}[theorem]{Fact}
\newtheorem{remark}[theorem]{Remark}

\newtheorem{observation}[theorem]{Observation}
\newtheorem{notation}[theorem]{Notation}

\newtheorem{conclusion}[theorem]{Conclusion}

\newtheorem{question}[theorem]{Question}
\newtheorem{stheorem}{Theorem}[subsection]

\newtheorem{slemma}[stheorem]{Lemma}
\newtheorem{sdefinition}[stheorem]{Definition}
\newtheorem{sclaim}[stheorem]{Claim}
\newtheorem{sconstruction}[stheorem]{Construction}
\newtheorem{sfact}[stheorem]{Fact}

\newtheorem{snotation}[stheorem]{Notation}

\def\b1K{\mbox{\boldmath $K$}_{-1}}
\def\bK{\mbox{\boldmath $K$}}

\def\Mscr{{\MM }}

\def\Ascr{{\mathcal {A} }}

\def\Iscr{{\mathcal {I} }}
\def\Pscr{{\mathcal {P} }}

\def\Tscr{{\mathcal {T} }}

\def\uf{\mathop{\rm uf}}

\def\meet{\wedge}

\def\bB{\mbox{\boldmath $B$}}

\newbox\noforkbox \newdimen\forklinewidth
\forklinewidth=0.3pt \setbox0\hbox{$\textstyle\smile$}
\setbox1\hbox to \wd0{\hfil\vrule width \forklinewidth depth-2pt
 height 10pt \hfil}
\wd1=0 cm \setbox\noforkbox\hbox{\lower 2pt\box1\lower
2pt\box0\relax}
\def\unionstick{\mathop{\copy\noforkbox}\limits}
\def\nonfork_#1{\unionstick_{\textstyle #1}}
\def\nonfork_#1{\unionstick_{\textstyle #1}}
\setbox0\hbox{$\textstyle\smile$} \setbox1\hbox to \wd0{\hfil{\sl
/\/}\hfil} \setbox2\hbox to \wd0{\hfil\vrule height 10pt depth
-2pt width
               \forklinewidth\hfil}
\newbox\doesforkbox
\setbox\doesforkbox\hbox{\lower 2pt\box1 \lower
2pt\box2\lower2pt\box0\relax}

\def\add{\widehat{\ }}

\newcommand{\BB}{\mbox{\msbm B}}
\newcommand{\FF}{\mbox{\msbm F}}

\newcommand{\MM}{\mbox{\msbm M}}

\def\sub'm{\prec_{\bK'}}
\def\grpf #1 #2{{\rm grp}_{#2}(#1)}
\def\spanf #1 #2{{\rm span}_{#2}(#1)}
\def\fldf #1 #2{{\rm fld}_{#2}(#1)}
\def\dclf #1 #2{{\rm dcl}_{#2}(#1)}
\def\rclf #1 #2{{\rm rcl}_{#2}(#1)}
\def\aclf #1 #2{{\rm acl}_{#2}(#1)}
\def\acff #1 #2{{\rm acf}_{#2}(#1)}
\def\strf #1 #2{{\rm str}_{#2}(#1)}
\def\tclf #1 #2{{\rm acf}_{#2}(#1)}
\def\abar{\mbox{\boldmath $a$}}
\def\bbar{{\bf b}}
\def\cbar{{\bf c}}

\def\hbar{{\bf h}}

\def\tbar{{\bf t}}

\def\wbar{{\bf w}}
\def\xbar{{\bf x}}
\def\ybar{{\bf y}}

\def\bT{\mbox{\boldmath $T$}}
\def\bF{\mbox{\boldmath $F$}}

\date{\today}
\begin{document}
\title{Maximal models up to the first measurable in  ZFC}

\author{John T. Baldwin
\\University of Illinois at Chicago\thanks{Research partially supported by Simons
travel grant G5402, G3535.}
\\
Saharon Shelah\thanks{Item 1147 on Shelah's publication list.
Partially supported by European Research Council grant 338821, and by National Science Foundation grant  136974 and DMS 1833363.
}\\Hebrew University of Jerusalem}

\maketitle

\begin{abstract}{ Theorem: There is a {\em complete sentence}
$\phi$ of $L_{\omega_1,\omega}$ such that $\phi$ has maximal models in a set of cardinals $\lambda$ that is cofinal
in the first measurable $\mu$ while $\phi$ has no maximal models in any
$\chi \geq \mu$.}
\end{abstract}

In this paper we prove in ZFC the existence of a {\em complete sentence}
$\phi$ of $L_{\omega_1,\omega}$ such that $\phi$ has maximal models in a set of cardinals $\lambda$ that is cofinal
in the first measurable $\mu$ while $\phi$ has no maximal models in any
$\chi \geq \mu$. In \cite{BaldwinShelahhanfmax}, we proved a theorem with a similar result;
the earlier proof required that $\lambda = \lambda^{< \lambda}$, and extended ZFC
 by requiring an $S \subseteq S^{\lambda}_{\aleph_0}$, that is stationary
non-reflecting, and $\diamond_S$ holds.  Here, we show {\em in ZFC} that
the sentence $\phi$ defined in \cite{BaldwinShelahhanfmax} has maximal models cofinally in  $\mu$.
 The additional hypotheses in
\cite{BaldwinShelahhanfmax} allow one to demand that if $N$ is a submodel with
cardinality $<\lambda$ of  the $P_0$-maximal model, $N$  is $\bK_1$-free
(See Remark~\ref{correnum2}); that property fails for the example here.
 The existence of such a $\phi$ which is {\em not complete} is well-known (e.g. \cite{Magcrm}).

 This paper contributes to the study of Hanf numbers for infinitary logics.
 Works such as   \cite{BKS,BKSoul,BSoul,KLH} study the spectrum of maximal
 models in the context where the class has a bounded number of models.
We list now some properties that are true in every cardinality for first order logic but
are true only eventually for complete sentences of $L_{\omega_1,\omega}$ or,
more generally, for abstract elementary classes, and compare the cardinalities (the Hanf number) at which the cofinal
 behavior must begin. Every infinite model of a first order theory
has a proper elementary extension and so each theory has arbitrarily large models.
Morley \cite{Morley65a} showed that every sentence of $L_{\omega_1,\omega}$ that
   has models up to $\beth_{\omega_1}$ has arbitrarily
    large models and provided counterexamples showing that cardinal was minimal.
Thus he showed the Hanf number for existence of  $L_{\omega_1,\omega}$-sentences
  in a countable vocabulary is
$\beth_{\omega_1}$. Hjorth \cite{Hjorthchar}, by a much more complicated argument,
 showed  there are {\em complete} sentences  $\phi_{\alpha}$ for $\alpha< \omega_1$ such
that $\phi_\alpha$ has a model in $\aleph_\alpha$ and no larger so the
Hanf number for complete sentences is $\aleph_{\omega_1}$.  The amalgamation
 property holds for every complete first order theory. However, \cite{BBHanf}
 show that an upper bound on the Hanf number
  for amalgamation is the first strongly compact; the actual value remains open.
Boney and Unger \cite{BoUn}, building on \cite{Sh932} show that the Hanf number
 `for all AEC's are tame' is the first strongly compact cardinal.
They also show the analogous property for various variants on tameness
 is equivalent to the existence of almost (weakly) compact, measurable, strongly compact).
 The result here shows in ZFC that the Hanf number for extendability
 (every model of a complete sentence has a proper $L_{\omega_1,\omega}$-elementary extension) is the first measurable cardinal.

Section~\ref{prelim} provides some background information on Boolean algebras. Section~\ref{stba} is a
set theoretic argument for the existence of a Boolean algebra with certain
specified properties in any cardinal
$\lambda$  of the form $\lambda = 2^\mu$ that is less than the first measurable;
 this construction is completely independent of the model theoretic results. Then we make the connection with model
 theory. In particular, we link the construction here with the complete sentence $\phi$ from \cite{BaldwinShelahhanfmax}.
 Section~\ref{1approx} builds several approximations to the counterexample.
 Subsection~\ref{fgmodels} introduces the most basic class  of models $\bK_{-1}$ and explains the
connections with \cite{BaldwinShelahhanfmax}.
 Subsection~\ref{mtconst} builds on this result to find a $P_0$-maximal model in $\bK_{-1}$ with
cardinality $\lambda$ satisfying certain further restrictions.
We recall in Subsection~\ref{K1def} the class $\bK_2$ of models of the complete sentence
 from \cite{BaldwinShelahhanfmax}.
In Section~\ref{corr}, the $P_0$-maximal model from Section~\ref{mtconst} is converted to the $P_0$-maximal
model in $\bK_2$. From this, it is easy to find a maximal model in $\bK_2$ of roughly the same cardinality.

The first author acknowledges helpful conversations with Joel Berman, Will Boney,
 Ioannis Souldatos, and especially Sherwood Hachtman.
We are particulary grateful for an extremely helpful referee report.

\section{Preliminaries}\label{prelim}

This paper depends heavily on \cite{BaldwinShelahhanfmax} which contains a fuller background and
 essential material on Boolean algebras. In particular, the incomplete sentence with maximal models cofinal
 in the first measurable is described there and the construction of the desired complete sentence; in this paper we show
 in ZFC that sentence has maximal models below the first measurable.
 We repeat in this section the main slightly nonstandard
definitions  from Boolean algebra that appear in \cite{BaldwinShelahhanfmax} and some immediate consequences.

\begin{definition}\label{bn}

\begin{enumerate}
\item A Boolean polynomial $p(v_0,\ldots ,v_k)$ is a term formed by the compositions of
the $\wedge, \vee, ^{-1}, 0,1$ on the variables $v_i$; a polynomial over $X$ arises when elements of $X$ are substituted for some of the $v_i$.
\item For $X \subseteq B$ and $B$ a Boolean algebra, $\overline X
    = X_B =\langle X\rangle_B$ denotes the subalgebra of $B$ generated
    by $X$.
\item A set $Y$ is {\em independent} (or {\em free}) over $X$
  {\em modulo} an ideal $\Iscr$ (with domain $I$) in a Boolean algebra $B$ if and only if
    for any Boolean polynomial $p(v_0,\ldots ,v_k)$ (that is not
    identically $0$, i.e. non-trivial), and any
   $a\in \langle X\rangle_B- \Iscr$, and distinct $y_i \in Y$,
 $p(y_0,\ldots ,y_k)\wedge a\not \in \Iscr$.

 \item A $Y$ which is independent over $X$  modulo $I$ is called  a {\em basis} for $\langle X \cup Y \cup I\rangle$ over $\langle X \cup I\rangle$.
\end{enumerate}
\end{definition}

In this context, `independent from' may sometimes be written `independent over'.  These are distinct notions for forking independence.


\begin{observation}\label{obs1a}

 If $\Iscr$ is the $0$ ideal, (i.e., $Y$  is independent over
    $X$),
\begin{enumerate}
\item the condition becomes:  for any $b\in \langle
    X\rangle_B- \{0\}$, $B \models p(y_0,\ldots ,y_k)\wedge b>
    0$. %
    That is, every finite Boolean combination of elements of
       $Y$ has non-empty meet with each non-zero $a \in \langle X\rangle_B$.

\item or, there is no non-trivial polynomial $q(\ybar,\xbar)$ and $\bbar \subseteq X$ such that
$q(\ybar,\bbar) =0$.

%
%
%
\end{enumerate}

That 2) implies 1) is obvious. For the converse, put a counterexample
$q(\ybar,\bbar) =0$ in disjunctive normal form.
 Then  for each disjunct (i.e. each constituent conjunction)
$q'(\ybar,\bbar) =0$ (not all variables of $q$ may appear in $q'$.) We can replace
those $b$'s that appear in $q'$ by a single element $b$  of $\langle X\rangle$
to get a $q''(\ybar,b)= 0$; $q''$ contradicts condition 1).
\end{observation}

With Observation~\ref{obs1a} we obtain an analog for Boolean algebras of
the notion of dependence in vector spaces in rings or fields: $\{y_0,\ldots ,y_k\}$
are dependent over $X$ if some non-trivial polynomial
$p(v_0,\ldots ,v_k, w_0, \ldots w_m)$ and some $\bbar$ from $X$, $p(\ybar,\bbar) =0$.
This yields that
if $B_2$ is freely generated over $B_1$, all atoms in $B_1$ remain atoms in $B_2$.  If
not, there would be an atom $a$ of $B_1$ and a term $\sigma(\bbar_2,\bbar_1)$
with $0_{B_1} < \sigma(\bbar_2,\bbar_1) < a$ and  $\sigma(\bbar_2,\bbar_1) \in B_1$.
 But then $B_2 \models \sigma(\bbar_2,\bbar_1) \wedge a =0$; this contradicts
the freeness assumption. This notion of dependence (a depends on $X$ if
and only if $a\in \langle X\rangle$) does not satisfy the exchange axiom.
 See \cite[Chapter 5]{Gratzer}
for the strong consequences if this dependence relation satisfies exchange.

There is no requirement that $\Iscr$ be contained in $X$.  Observe the following:

\begin{observation}\label{obs1} Let $\Iscr$ be an ideal  in a Boolean algebra $B$.
\begin{enumerate}

\item Let $\pi$ map $B$ to $B/\Iscr$. If `Y is independent from
    $X$ over $\Iscr$' then the image of $Y$ is free from the
    image of $X$ (over $\emptyset$) in $B/\Iscr$.  Conversely, if
    $\pi(Y)$ is independent over $\pi(X)$ in $B/\Iscr$, for any
    $Y'$ mapping by $\pi$ to $\pi(Y)$, $Y'$ is independent from
    $X$ over $\Iscr$.

So, if  $X$ is empty, the condition `Y is independent  over
    $\Iscr$' implies  the image of $Y$ is an independent subset
    of $B/\Iscr$.

    \item If a set $Y$ is {\em independent} (or {\em free}) from
        $X$ over  $\Iscr$ in   $B$ and
        $Y_0$ is a subset of $Y$, then $Y-Y_0$ is {\em
        independent} (or {\em free}) from $X \cup Y_0$ ($\langle
        X \cup Y_0 \rangle_B$) over the ideal $\Iscr$ in the
        Boolean algebra $B$.


\end{enumerate}
\end{observation}

\section{Set theoretic construction of a Boolean algebra}\label{stba}

We define a property $\boxplus(\lambda)$, which asserts the existence in $\lambda$ of a Boolean algebra
that is `uniformly $\aleph_1$-incomplete'. We then show certain conditions on $\lambda$ imply $\boxplus(\lambda)$.
So this section has no model theory.
The arguments here are similar to those around page 7 of \cite{GobelShelah}.
 We connect this construction with
our model theoretic approach in
Section~\ref{1approx}.

\begin{definition}[$\boxplus(\lambda)$] denotes: \label{boxplus}
There are a Boolean algebra $\BB \subset \Pscr(\lambda)$ with $|\BB| = \lambda$ and a set $\Ascr \subseteq {}^{\omega} \BB$ such that:
\begin{enumerate}[i)] \item $\Ascr$  has cardinality $\lambda$ and if
$\overline{A} =\{A_n: n \in \omega\} \in \Ascr$ then for $\alpha < \lambda$ for all but finitely many $n$, $\alpha \not \in A_n$.
\item  $\BB$
     includes the finite subsets of $\lambda$; but is such that for every
      non-principal ultrafilter $D$ of $\lambda$ (equivalently an ultrafilter of $\BB$ that is
      disjoint from $\lambda^{<\omega}$) for some sequence $\langle A_n :
      n \in \omega\rangle \in \Ascr$, there are infinitely many $n$ with $A_n \in D$.

    \end{enumerate}
    \end{definition}

We may say that $(\BB, \Ascr)$ witness uniform $\aleph_1$-incompleteness.


\begin{theorem}[ZFC]\label{boxthm}
      Assume  for some $\mu$, $\lambda = 2^\mu$ and $\lambda$ is less than the
       {\em first measurable}, then  $\boxplus(\lambda)$ from \ref{boxplus} holds.
      \end{theorem}

We need the following structure.

      \begin{definition}\label{f12.5}
        \begin{enumerate}
        \item Fix the vocabulary
        $\tau$ with unary predicates $P,U$, a binary predicate $C$, and a binary function $F$.

      \item    Let $\langle C_\alpha\mcolon \alpha<\lambda\rangle$ list {\em without
      repetitions} $\Pscr(\mu)$ such that $C_0 = \emptyset$ and also let $\langle f_\alpha\mcolon
      \mu \leq \alpha<\lambda\rangle$ list\  ${}^\mu \omega$.
      \item Define the $\tau$-structure $M$
       by:
\begin{enumerate}
\item The universe of $M$ is $\lambda$;
\item $P^M = \omega$; $U^M = \mu$;

 \item  $C(x,y)$ is binary relation
  on $U \times M$ defined by $C(x,\alpha)$ if and only $x \in C_\alpha$.
  Note that $C$
 is extensional. I.e., elements of $M$ uniquely code subsets of $U^M$;
\item Let $F^M_2(\alpha, \beta)$ map $M \times U^M \rightarrow P^M$ by $ F^M_2(\alpha, \beta)  =f_{\alpha}( \beta)$ for $\alpha <\lambda$, $\beta< \mu$;
 \item $F^M_2(\alpha, \beta) = 0$ for $\alpha < \lambda$ and $\beta \in [\mu,\lambda)$.
     \end{enumerate}
      \end{enumerate}
     \end{definition}



   We  use the following, likely well-known, fact pointed out to us by Sherwood Hachtman.

   \begin{fact}\label{hacht} Let $D \subseteq \Pscr(X)$ and suppose that for each partition $Y\subseteq \Pscr(X)$ of $X$ into at most countably many sets, $|D\cap Y| = 1$.
 Then, $D$ is a countably complete ultrafilter.
 \end{fact}




   We use the following lemma about $M$ to find a Boolean algebra $\BB$ in $M$ that
   satisfies $\boxplus$. We lay the basis for the notion of $P$-maximality, a counterexample to
   maximality must occur in a given predicate $P$ (Definition~\ref{f2*}).

 \begin{lemma}\label{f12.7} If $\lambda$ is less than the first measurable cardinal and $\lambda = 2^{\mu}$ for some $\mu$ there is a model $M$, with  $|M| =\lambda$,  and a countable vocabulary with  $P^M$ denoting the natural numbers  such that  every first order proper elementary extension $N$ of $M$ properly extends $P^M$.
\end{lemma}

Proof. Fix  $M$  as in Definition~\ref{f12.5}. 
 We first show that any proper elementary extension $N$ of $M$ extends $U^M$.
Suppose for contradiction there exists $\alpha' \in N-M$ but $U^N =U^M$. By
the full listing of the $C_\alpha$, there is a   $\beta \in  M$ with $\{x: N \models C(x,\beta)\} = \{x: N \models C(x,\alpha')\}$. This contradicts extensionality of the relation $C$ in $N$; but $C$ is extensional in the elementary submodel $M$.

Now we show that if $U^M \subsetneq U^N$ and $P^M = P^N$, then there is a countably complete
non-principal ultrafilter on $\mu$, contradicting that $\mu$ is not measurable.
    Note that the sequence $\langle f_\alpha\colon \mu \leq \alpha< \lambda\rangle$ can be viewed as a list of all  non-trivial partitions of $\mu$ into at most countably many pieces.  Let $\nu^* \in U^N-U^M$. For $\alpha \in N$, denote 
   $F^N_2(\alpha,\nu^*)$ by $n_\alpha$. Since $P^M = P^N$, $n_\alpha\in M$.
   By elementarity, for $\alpha\in M,\eta \in U^M$, $F^N_2(\alpha,\eta) =
     F^M_2(\alpha,\eta)= f_\alpha(\eta)$.
    Now, let 


   $$D = \{x \subseteq U^M\colon  x \neq \emptyset \wedge (\exists \alpha \in M)
    \  x \supseteq f^{-1}_\alpha(n_\alpha)\}.$$

We show $D$ satisfies the conditions from Fact~\ref{hacht}.
 Let $W$ be a partition, indexed by $f_\alpha$.
 Then $f^{-1}_\alpha(n_\alpha) \neq \emptyset$ and is in $D$. Suppose for
contradiction there are $x_0\neq x_1$ in $W$ that are both in $D$.
 Then, there are $\alpha_i \in M$ such that $x_i \in W\cap D$ contains $f^{-1}_{\alpha_i}(n_{\alpha_i})$ for $i= 0,1$.
   So,
$N \models F(\alpha_i,\nu^*) = n_{\alpha_i}$ for $i =1,2$.  Since $\alpha_i \in M$ and $M \prec N$, $M\models \exists x (F(\alpha_0,x) = n_{\alpha_0} \wedge F(\alpha_1,x) = n_{\alpha_1}$.  So, by Definition~\ref{f12.5} (d),
for any witness $a$ in $M$ for this formula, $a \in x_0\cap x_1$;  but $x_0\cap x_1 =\emptyset$ since $W$ is a partition.

  Finally, $D$ is non-principal on $U^M$ since if it were generated by an $a\in U^M$,  $$D = \{x \subseteq U: (\exists \alpha)\   x \supseteq f^{-1}_\alpha(n_\alpha)\} =
   \{x \subseteq U:a \in x\}.$$
     Since $\{a\} \in D$, for some $\alpha_0 \in M$, $\{a\} = f^{-1}_{\alpha_0}(n_{\alpha_0})$. Note that $\alpha_0 \in M
   $, because the definition of $D$ is about the model $M$. That is,
   $M\models \exists ! y F(\alpha_0, y) = n_{\alpha_0}$.
    But $N \models F(\alpha_0, a) = n_{\alpha_0} \wedge F(\alpha_0, \nu^*) = n_{\alpha_0}$. This contradicts the assumption $M \prec N$ and completes the proof.
    $\qed_{\ref{f12.7}}$

\medskip

The following claim completes the proof of Theorem~\ref{boxthm}

 \begin{claim}\label{f12.8} If $\BB$ is the Boolean algebra of definable formulas
  in the $M$ defined in Definition~\ref{f12.5}, there is an $\Ascr$ such that
  $(\BB,\Ascr)$ is uniformly $\aleph_1$-incomplete so  $\boxplus(\lambda)$ holds.
\end{claim}
  Proof. We may assume $\tau$ has Skolem functions for $M$  and then define $\BB$
   and $\Ascr$ as follows to satisfy $\boxplus$.(ii).  Let $\BB$ be the Boolean
   algebra of definable subsets of $M$. I.e.,
    $$\BB =\{X \subseteq M\mcolon \ \textrm{\ for some }\ \tau\textrm{-formula}\
     \phi(\xbar,\ybar) \textrm{ and}\  \bbar \in {}^{\lg(\ybar)}M,\  
     \phi(M,\bbar) =X.\}$$

     Note $\BB$ is a Boolean algebra of cardinality $\lambda$ with the normal operations.
We define the Skolem functions a little differently than usual: as  maps
 $\sigma_\phi = \sigma_{\phi(x,w,\ybar)}$ from $M^{n+1}$ to $M$ for  formulas $\phi(x,w,\ybar)$ such
 that $\phi(\sigma_{\phi}(b,\abar),b,\abar)$. Here $\lg(\ybar) = n$. Then, we
 specialize the Skolem functions by considering the unary function arising
 from fixing the $\ybar$ entry of  $\sigma_{\phi}(w,\ybar)$ to obtain
 $\sigma_{\phi}(w,\abar)$.

%



\begin{eqnarray*}
A_n^{\sigma_{\phi}(w,\abar)} & = &\{\alpha< \lambda \mcolon: \phi(\sigma^M_\phi(\alpha,\abar),\alpha,\abar) \wedge 
P(\sigma^M_\phi(\alpha,\abar)) \wedge \sigma^M_\phi(\alpha,\abar) \nleqslant n\}\\
 & \cup & \{\alpha< \lambda \mcolon:\ n=0 \wedge \neg P(\sigma^M_\phi(\alpha,\abar)\}.
\end{eqnarray*}
%
%

Then let
   $ \overline{A}_{\sigma_\phi(w,\abar)} =\langle A_n^{\sigma_\phi(w,\abar)}\mcolon n<\omega\rangle$ and
 $$ (*) \hskip .1in
 \Ascr =\{\overline{A}_{\sigma_\phi(w,\abar)}: \ \textrm{ for some} \   \tau_M\textrm{-term}\  \sigma_\phi(w,\ybar)\
\textrm{ and} \ \abar \in {}^{\lg(\ybar)}M. \}$$
Note $|\Ascr| = \lambda =\lambda^\omega$ as for each $\abar \in M$ and each of the countably many terms $\sigma_\phi(w,\abar)$, $\overline{A}_{\sigma_\phi(x,w,\abar)}$ is  a map from $\omega$ into $\BB$.
For each $\alpha$,
for each $0<m < \omega$ and $\overline{A}= \overline{A}_{\sigma_\phi(\alpha,\bbar)}$, the set
 $\{m \colon \alpha \in A_m\}$ is finite, bounded by $\sigma_\phi(\alpha,\abar)$.
  Thus, clause i) of $\boxplus$ is satisfied.

 We now show Clause ii) of $\boxplus$.  Let $D$ be an arbitrary non-principal
 ultrafilter on $\lambda$ and  where
 $ \phi(v,\ybar)$ varies over first order $ \tau$-formulas such that $\ybar$ and $\abar$ have
 the same length, define the type $p(x)=p_D(x)$ as:
$$p(x) =\{\phi(x,\abar)\colon \{\alpha\in M\colon  M \models \phi(\alpha,\abar)\} \in D\}.$$
 Since $D$ is an ultrafilter, $p$ is a complete type over
 $M$.  So there is an elementary extension $N$ of $M$ where an element
 $d$ realizes $p$.   Let $N$ be the Skolem hull of $M \cup \{d\}$.
 Since $D$ is non-principal, so is $p$; thus, $N \neq M$. By Lemma~\ref{f12.7}, we can
 choose a witness $c \in   P^N -P^M $.
 Since, $N$ is the Skolem hull of $M \cup \{d\}$ there is a Skolem term $\sigma(w,\ybar) = \sigma_\phi(w,\ybar) $ and $\abar \in M$ such that $c = \sigma^N(d,\abar)$.
 Since $c \not \in M$,  for each $n\in P^M$,
 $N      \models \bigwedge_{k<n} c \neq k$ so
 $N                  \models \bigwedge_{k<n} \sigma(d,\abar) \neq k$
 so
 $\bigwedge_{k<n} \sigma(x,\abar) \neq k$ is in $p$. That is, for each $\sigma_\phi$ and each $n$, $A^{\sigma_\phi(w,\abar)}_n$
  is in $D$.
  $\qed_{\ref{f12.8}}$

\section{Three Classes of Models and an Approximate Counterexample}\label{1approx}
\numberwithin{theorem}{subsection}
In this section we define the model theoretic classes that produce
first an amalgamation class $\bK_{-1}$ of finite structures (Section~\ref{fgmodels}),
then the class $\bK_2$  (Definition~\ref{richname}) of models of  a complete
$L_{\omega_1,\omega}$-sentence.  Using Theorem~\ref{boxthm}, we build in Subsection~\ref{mtconst} a model $M_*$
 in  $\bK_{-1}$ with cardinality $\lambda$, which is $P_0$-maximal.
 Subsection~\ref{K1def} defines the classes $\bK_1$ and $\bK_2$ which give us the complete sentence.
In Section~\ref{corr} we modify $M_*$ to a  $P_0$-maximal model in $\bK_2$ and then construct the required
maximal model in  $\bK_2$.

\subsection{Finitely generated models}\label{fgmodels}

We include, as needed, definitions of the classes of model $\bK_{-1},\bK_{1}, \bK_{2}$
introduced in \cite{BaldwinShelahhanfmax}.  For each $i$, $\bK^{i}_{<\aleph_0}$ denotes the
class of finitely generated  members of $\bK_i$.

\begin{definition}\label{deftau} $\tau$ is a vocabulary with unary predicates $P_0, P_1,
    P_2, P_4$, binary $R$,  $\wedge, \vee$, $\leq$ 
     unary
    functions $^-$, $G_1$,
    constants 0,1 and unary functions $F_n$, for $n< \omega$.  $\leq$ is a partial order on $P^M_1$ and the Boolean algebra can be defined from it.
\end{definition}

We occasionally use the notations $(\forall^{\infty} n)$ and $(\exists^{\infty} n)$
to mean `for all but finitely many' and `for infinitely many' respectively.
  It is easy to see that $\bK_{-1}$ is $L_{\omega_1,\omega}$-axiomatizable but far from complete.

\begin{definition}[$\bK_{-1}$]\label{f1}
$M \in \bK^{-1}_{< \aleph_0}$ is the class of {\em finitely generated} structures $M$ satisfying the following conditions.


\begin{enumerate}
\item $P^M_0, P^M_1, P^M_2$ partition $M$.
\item $(P^M_1, 0,1, \meet, \vee,\leq,^-)$ is a Boolean algebra
    ($^-$ is complement).
    We also
   consider ideals and restrictions to them of the relations/operations except for complement.

\item $R \subset P^M_0 \times P^M_1$ with $R(M,b)=
    \{a:R^M(a,b)\}$  and the set of $\{R(M,b): b \in P^M_1\}$ is
    a Boolean algebra. $f^M\colon P^M_1\mapsto \Pscr(P^M_0)$ by
    $f^M(b) =R(M,b)$ is a Boolean algebra homomorphism into
    $\Pscr(P^M_0)$.

    Note that $f$ is not\footnote{The subsets of $P^M_0$ are {\em
    not} elements of $M$.} in $\tau$; it is simply a convenient
    abbreviation for the relation between the Boolean algebra
    $P_1^M$ and the set algebra on $P_0$ by the map $b \mapsto
    R(M,b)$.

    \item $P^M_{4,n}$ is the set containing each join of $n$
        distinct atoms from $P^M_1$; $P^M_{4}$ is the union of the
        $P^M_{4,n}$ and so is an ideal. That is, $P^M_4$ is the set of all finite
        joins of atoms.

   There is an element     $b^* \in P^M_1$ such that $P^M_{4} = \{c:c \leq^M b_*\}$.
Note that $b_*$ is not a function symbol in $\tau$.



%


\item $G^M_1$ is a bijection from $P^M_0$ onto $P^M_{4,1}$ such
    that $R(M,G^M_1(a)) = \{a\}$. (Note that $P^M_0 = \emptyset$ is allowed.


    \item $P^M_2$ is finite
   (and may be empty). Further, for each $c \in
      P^M_2$ the $F^M_n(c)$ are functions from  $P^M_2$ into $P^M_1$.
      Note that it is allowed that for all but finitely many $n$, $F^M_n(c) = 0_{P^M_1}$.
%
%

\item ({\em countable incompleteness}) If $a \in P^M_{4,1}$ and $c\in P^M_2$ then
 $(\forall^{\infty} n)\  a\nleqslant_M F^M_n(c)$. Since $a \wedge F^M_n(c) =0$ and
 $a$ is an atom, this implies  $\bigwedge_{n\in \omega}
\{x\colon(G_1(x)  \in F^M_n(c)\} = 0$.

         \item $P^M_1$ is generated
as a Boolean algebra by $P^M_4 \cup \{F^M_n(c)\mcolon c \in P^M_2, n\in \omega\} \cup X$ where
$X$ is a finite subset of $P^M_1$.

%

    \end{enumerate}


\end{definition}

\begin{definition}\label{f2}
\begin{enumerate}
\item $K_{-1}$ is the class of $\tau$ structures $M$ such that every
finitely generated substructure of $M$ is in   $\bK^{-1}_{< \aleph_0}$.
$K^{-1}_\mu$ is the members of $K_{-1}$ with cardinality $\mu$.
    \item We say $M \in K_{-1}$ is {\em atomic} if $P^M_1$ is
    atomic as a Boolean algebra. That is, $P^M_4$ is dense in $\bB_M$.
 \end{enumerate}
 \end{definition}

\bigskip
\subsection{A $P_0$-maximal model in $\bK_{-1}$}\label{mtconst}

In this section we invoke Theorem~\ref{boxthm} to show (Theorem~\ref{f11a}) that we can construct $P_0$-maximal structures
in the class $\bK_{-1}$ of appropriate cardinality below the first measurable.


\begin{sdefinition}\label{f2*} We say $M \in K_{-1}$ is {\em $P_0$-maximal} (in $K_{-1}$)
         if $M \subseteq N$ and $N \in K_{-1}$ implies $P^M_0 =P^N_0$.

            \end{sdefinition}


The notion $\uf(M)$ is the crucial link between Section~\ref{stba} and $P_0$-maximality.
Lemma~\ref{f8} is central for Theorem~\ref{f11a} and is applied in Theorem~\ref{f44}.

\begin{sdefinition}[$\uf(M)$]\label{f5} For $M\in K_{-1}$, let
${\rm uf}(M)$ be the set of  ultrafilters $D$ of the Boolean Algebra $P^M_1$
  such that $D \cap P^M_{4,1} = \emptyset$ and for each $c \in P^M_2$
  only finitely many of the $F^M_n(c)$ are in $D$.

For applications we rephrase  this notion with the following  terminology.
For any $M \in \bK_{-1}$ and $d \in P^M_2$, let $S^M_d(D) =\{n: F^M_n(d) \in D\}$. So $\uf(M) = \emptyset$ if and only if for every ultrafilter $D$ on $P^M_1$,
there exists a $d \in P^M_2$ such that  $S^M_d(D)$ is infinite.

\end{sdefinition}

%

We use the following standard properties of a  Boolean algebra $B$ and ideal
$I$  in proving Lemma~\ref{f8} and
deducing  Claim~\ref{succ} from Definition~\ref{defI}.

\begin{sfact}\label{quotprop}
\begin{enumerate}
\item $b\wedge c \in I$ implies $b/I$ and $c/I$ are disjoint.
\item $b \vartriangle c \in I$ implies $b/I =c/I$.
\item $b -c \in I$ implies $b/I \leq c/I$.
\end{enumerate}
\end{sfact}

For our collection of structures $\bK_{-1}$, we can characterize {$P_0$}-maximality in terms of ultrafilters.


       \begin{slemma}\label{f8}  An $M \in \bK_{-1}$ is $P_0$-maximal if and
       only if $\uf(M) =\emptyset$.
       \end{slemma}

       Proof.  Suppose $M$ is not $P_0$-maximal and $M\subset N$ with $N \in \bK_{-1}$
       and $d^* \in P_0^N - P_0^M$.
        Then $\{b \in M: R^N(d^*,b)\}$ is a non-principal ultrafilter
       $D_0$ of the Boolean algebra $P^{M}_1$ \cite[3.3.11]{BaldwinShelahhanfmax}.
        To see $D_0$ is non-principal
       suppose there is a $b_0 \in P^M_1$ such that $D_0 = \{b\in M: b_0 \leq b\}$.
       Note $b_0 = G^M_1(a)$ for some $a \in P^M_0$.
       But $N \models G^N_1(d^*) \ngeq b_0$, contradicting $\{d^*\} \in D_0$.

        For each $c \in P^M_2$, since $N\in \bK_{-1}$, by countable incompleteness
         (clause 7 of Definition~\ref{f1}),
       for all $a \in P^N_0$ and all but finitely many $n$, $G^N_1(a) \not \leq F^N_n(c)$.
        Since $F^N_n(c) = F^M_n(c)$,
          only finitely many of the $F^M_n(c)$ can be in $D_0$, which implies $D_0 \in \uf(M)$. By contraposition we have the right to left.

Conversely, if $D \in \uf(M)$, we can construct an extension
 by adding an element $d \in  P_0^N$ satisfying $R^N(d,b)$ iff $b \in D$.
Let $P^N_1$ be the Boolean algebra generated by $P^M_1\cup \{G_1(d)\}$
modulo the ideal generated by $\{G^N_1(d)-b :b\in D\}$;
 this implies that in the quotient $G_1(d)\leq b$.
  (Compare Fact~\ref{quotprop}).
  Let $P_2^N= P_2^M$ and $F^N_n(c) = F^M_n(c)$.
   Since $D \in \uf(M)$, it is easy to check that $N \in \bK_{-1}$.
 $\qed_{\ref{f8}}$

 \medskip
 We now introduce the requirement that the Boolean algebras constructed will,
  when the atoms are factored out, be free. Moreover,
  there is a set $Y\subseteq P^N_2$ with $|Y|=\lambda$ such that different $c\in Y$
generate coinitially disjoint collections of $F^N_n(c)$ as $c$ varies. This strong requirement
is used inductively in this section to construct an approximation to the counterexample.
The correction in Section~\ref{corr}
loses this disjointness (and thus freeness).

       \begin{sdefinition}[Nicely Free]\label{b9}
       We say $M \in K_{-1}$ is {\em nicely free} when $|P^M_1| = \lambda$ and there is a sequence $\bbar = \langle b_\alpha \mcolon \alpha< \lambda \rangle$ such that
       \begin{enumerate}[(a)]
       \item $b_\alpha \in P^M_1 - P^M_4$;
       \item $\langle b_\alpha/P^M_4 \mcolon \alpha< \lambda
       \rangle$ 
           generate $P^M_1/P^M_4$ freely;
           \item 
there is  a set $Y\subset P^M_2$ of cardinality $\lambda$ such that
$\{F_n(c) : n < \omega; c \in Y\}$ without repetition is
a subset of the basis $\{b_\alpha : \alpha < \lambda\}$ mod atoms.  For $c\in Y$,
we  write
$u_c = \{F^M_n(c)\colon n < \omega\}$.

%
\end{enumerate}
\end{sdefinition}

Nicely free is quite distinct from the notion {\em $\bK_1$-free} introduced in
\cite{BaldwinShelahhanfmax}.  There are maximal nicely free models but there are no
maximal $\bK_1$-free models. Note that condition Definition~\ref{b9}.c asserts that a {\em sub}set
of $P^{M}_2$ partitions a {\em sub}set of the basis.

\medskip

  Here is the main theorem of Section~\ref{1approx}.
The hypotheses $\lambda =2^\mu$ and $\lambda$ is less than the first measurable cardinal  were used essentially
  as  the hypotheses for proving $\boxplus(\lambda)$, the existence of a uniformly $\aleph_1$-incomplete Boolean algebra.
   But here we use $\boxplus(\lambda)$ and don't rely again  on $\lambda$ being less than the first measurable cardinal.
  The  argument here does depends on $\lambda = \lambda^{\aleph_0}$, which follows from  $\lambda =2^\mu$.
   By constructing a nicely free model, we introduce
    at this stage the independence requirements,
    needed in Section~\ref{corr} to satisfy Definition~\ref{k0}.6, on the $F_n(c)$.


       \begin{stheorem}\label{f11a}
       If for some $\mu$, $\lambda =2^\mu$ and $\lambda$ is less than the
       first measurable cardinal  then
there is a  $P_0$-maximal model $M_*$ in $\bK_{-1}$ such that
$|P^{M_*}_i| = \lambda$ (for $i = 0,1,2$), $P^{M_*}_1$ is an atomic Boolean algebra,
$\uf(M_*) = \emptyset$, 
and $M_*$ is nicely free.
\end{stheorem}

Proof.
%
%
%
    We first construct by induction a $P_0$-maximal model in $\bK^{-1}$.    The
     property  $\boxplus(\lambda)$ (Definition~\ref{boxplus}) appears in
      the construction to satisfy Specification (f)
      and is used in the proof that the construction works in considering possibility
       2.
 We  choose $M_\epsilon, D_\epsilon$ and other auxiliaries by induction for $\epsilon \leq \omega+1$ to satisfy the following {\em specifications} of the construction.



   \begin{sconstruction}[Specifications]\label{oplus}
    \begin{enumerate}[(a)]
    \item For $\epsilon \leq \omega+1$, $M_\epsilon$ is a continuous increasing chain of
     members of $K^{-1}_\lambda$ with each $P^{M_\epsilon}_1$ atomic and $P^{M_{\omega+1}}_1= P^{M_{\omega}}_1 $.
     \item \label{clb} For all $\epsilon \leq \omega$, $|P^{M_\epsilon}_i| =\lambda$
         and $P^{M_{\omega}}_{i} =P^{M_{\omega+1}}_{i}$ for $i =0,1$.
         \item \label{clc} For all $\epsilon \leq \omega+1$, $P^{M_\epsilon}_1/P^{M_\epsilon}_4$ is a free Boolean algebra;
        \item \label{cld}{} \begin{enumerate}[(i)]
        \item If $\epsilon < \omega$, $D_\epsilon \in \uf(M_\epsilon)$.

  \item If $\epsilon =0$, then $\bbar_{-1} =\langle b_{-1,\alpha}\colon  \alpha < \lambda\rangle$ is a free basis of $P^{M_0}_1/P^{M_0}_4$, listed without repetition, and
                $\langle F^{M_0}_n(c)\colon n< \omega, c \in P^{M_0}_2\rangle$ lists $\langle b_{-1,\alpha}\colon  \alpha < \lambda\rangle$ without repetition.

        \item if $\epsilon = \zeta+1< \omega$ then there is a free basis
         $\bbar_{\zeta} =\langle b_{\zeta, \alpha}/P^{M_\zeta}_4 \colon \alpha< \lambda\rangle$
          of $P^{M_\epsilon}_1/P^{M_\epsilon}_4$.
          Note $b_{\zeta, \alpha} \in P^{M_\epsilon}_1 - P^{M_\zeta}_1$.

            \end{enumerate}

%

         \item \label{cle} if $\epsilon = \omega+1$, for each $\overline{d} \in
         {}^{\omega}(P^{M_{\omega+1}}_1-P^{M_{\omega+1}}_4)$ such that for each
          $a \in P^{M_\omega}_0$ satisfying that all but finitely many $n$,
          $a \not \in  R(M_\omega,d_n)$, there is a $c\in P^{M_{\omega+1}}_2$,
          $F^{M_{\omega+1}}_n(c) = d_n$;
(We will in fact have that $P^{M_{\omega+1}}_1 = P^{M_{\omega}}_1$ and
$P^{M_{\omega+1}}_4 = P^{M_{\omega}}_4$.)
%


             \item \label{clf} $\epsilon = \zeta+1 <  \omega$:

             Let $\BB$ and $\Ascr$ be as in Definition~\ref{boxplus}.
              There is a $1$-$1$ function $f_\epsilon$ from $\lambda$ onto $P^{M_{\epsilon}}_{4,1}$
               such that:
               \begin{enumerate}[i)] \item for every $X \in \BB$ (from $\boxplus$)
                there is a $b  = b_X \in P^{M_\epsilon}_1$ such that
                $$\{ \alpha< \lambda\mcolon f_\epsilon(\alpha) \leq_{M_\epsilon} b_X\} = X;$$

                 \item  
                     for each $\overline{A} =\langle A_n: n<\omega\rangle \in \Ascr$
                      there is a $c \in P^{M_\epsilon}_2$ such that for each $n$:
                     $$A_n=\{\alpha<\lambda\mcolon f_\epsilon(\alpha) \leq_{P^{M_\epsilon}_1} F^{M_\epsilon}_n(c)\}.$$

%

        \end{enumerate}

        \end{enumerate}
        \end{sconstruction}


\medskip
{\bf Carrying out the construction.}

\medskip

\underline{case 1:} When $\epsilon = 0$,
take $P^{M_0}_1$ as the Boolean algebra generated  by a set $P^{M_0}_{4,1}$ of cardinality $\lambda$ along with a set $ \{b_{-1,\alpha}\colon \alpha < \lambda\}$ of independent subsets of $\Pscr(\lambda)$. Let $G_1$ be a bijection between a set
$P^{M_0}_{0}$ and $P^{M_0}_{4,1}$. Set $P^{M_0}_{4}$ as the ideal generated by the image of $G_1$. For $a\in P^{M_0}_{0}$ and $b \in P^M_1$, define $R^{M_0}(a,b)$ to hold if $G_
1(a) \leq b$.  Set $P^{M_0}_2$ as a set of cardinality of $\lambda$ and let $\langle F^{M_0}_n(c)\colon n< \omega, c \in P^{M_0}_2\rangle$ list  $\langle b_{-1,\alpha}\colon  \alpha < \lambda\rangle$ without repetition.
Thus, any non-principal ultrafilter on $P^{M_0}_1$ is in
$\uf(M_0)$.

\underline{case 2:} For $\epsilon = \omega$, $M_\omega = \bigcup_{n<\omega} M_n$.
Since the set of free generators is extended at the step, the union is also free mod $P^M_4$.

\underline{case 3:} If $\epsilon = \zeta +1 < \omega$, the main effort is to
 verify clauses \ref{clc}, \ref{cld}, and \ref{clf} of Specification~\ref{oplus}.
The element $b_{\zeta,a_\alpha}$ is the $b_{A_\alpha}$ from Specification~\ref{oplus}.f.(i).

\smallskip
Now, to construct $M_\epsilon$:
\begin{enumerate}[(i)]
 \item Recall that $D_\zeta \in \uf (M_\zeta)$.
 \item   choose  as the new atoms introduced at this stage a set
 $B_\epsilon \subseteq \Pscr(\lambda)$
 with $B_\epsilon \cap M_\zeta = \emptyset$ and $|B_\epsilon| = \lambda$.

\item Let $f_\epsilon$ be a one-to-one function from $\lambda$ onto $B_\epsilon\cup P^{M_\zeta}_{4,1}$.

\item Let $\langle X_\gamma \mcolon \gamma < \lambda\rangle$ list the elements of $\BB$ (definable subsets of $M$ \ref{f12.8}) from $\boxplus$.(ii) with $X_0 = \emptyset$.

\item Fix a sequence  $\{b_{\zeta,\alpha}\mcolon \alpha < \lambda\}$, which are distinct and not in $M_\zeta \cup B_\epsilon$,  and let $\BB'_\zeta$ be the Boolean Algebra generated freely by $$P^{M_\zeta}_1 \cup \{b_{\zeta,\alpha}\mcolon \alpha < \lambda\} \cup \{f_\epsilon(\alpha): \alpha < \lambda\}.$$ 
\end{enumerate}

Using Lemma~\ref{quotprop}, we apply the following definition at the successor stage.
 Here we take an abstract Boolean algebra
 $\BB'_\zeta$ and impose relations  to embed $P^{M_\zeta}_1$ in a quotient
  $\BB''_\zeta$ of $\BB'_\zeta$.

\begin{sdefinition}[Ideal]\label{defI} Let $I_\zeta$ be the ideal of $\BB'_\zeta$ generated by:
\begin{enumerate}[(i)]
\item $\sigma(a_0, \ldots a_m)$ when $\sigma(x_0, \ldots x_m)$ is a Boolean
 term, $a_0, \ldots a_m \in P^{M_\zeta}_1$ and $P^{M_\zeta}_1 \models \sigma(a_0, \ldots a_m)= 0$.

    The next two clauses aim to show that in $M_\zeta/I_\zeta$,
    the element $b_{\zeta,\gamma}$ is the $b_{X_\gamma}$ from
    Specification~\ref{oplus} f.i).  That is,
    $\{ \alpha< \lambda\mcolon f_\epsilon(\alpha) \leq_{M_\epsilon} b_{\gamma,\zeta} \} = X_\gamma$.
 Recall (Definition~\ref{boxplus}) that the $X_\gamma$ enumerate $\BB$ and are subsets of $\lambda$.

    \item 

    $f_\epsilon(\alpha)- b_{\zeta,\gamma}$
  when $\alpha \in X_\gamma$
    and $\alpha,\gamma < \lambda$.

     \item   $b_{\zeta,\gamma} \wedge f_\epsilon(\alpha)$ when $\alpha \in \lambda - X_\gamma$
            and $\alpha,\gamma < \lambda$.


To show the $f_\epsilon(\gamma)$ are disjoint   atoms we add:

\item For any $f_\epsilon(\gamma)$ and any $b \in \BB'_\zeta$
either $(f_\epsilon(\gamma) \wedge b) \in I_\zeta$ or $(f_\epsilon(\gamma) - b)\in I_\zeta$.
         \item $f_\epsilon(\gamma_1) \wedge f_\epsilon(\gamma_2)$ when $\gamma_1 <\gamma_2 < \lambda$;

         \item $f_\epsilon(\alpha)- b$ when $\alpha < \lambda$,
          $f_\epsilon(\alpha) \not \in P^{M_\zeta}_{4,1}$ and $b \in D_\zeta$.

This asserts:             Every new atom is below each $b \in D_\zeta$ and is used at the end of case 3 of the construction.


    \end{enumerate}
\end{sdefinition}

Let $\BB''_\zeta =  \BB'_\zeta/I_\zeta$.  Applying Fact~\ref{quotprop}, we see from Definition~\ref{defI}:

\begin{sclaim}\label{succ} The structure $P^{M_\zeta}_1$ is embedded as a Boolean algebra into $\BB''_\zeta$ by the map $b \mapsto b/I_{\zeta}$ and
\begin{enumerate}
\item  For $\gamma < \lambda$, $f_\zeta(\gamma)/I_{\zeta}$  is an atom of $\BB''_\zeta$;
    \item  If $b \in P^{M_\zeta}_1$ is non-zero, then $b/I_{\zeta} \geq_{\BB''_\zeta} f_\epsilon(\gamma)$ for   some $\gamma< \lambda$.
    (Since $f_\epsilon^{-1}$ induces an isomorphism of ${\BB''_\zeta}$ into $\Pscr(\lambda)$

        \end{enumerate}

        \end{sclaim}

We take a further quotient of $\BB'_\zeta$.
       Let $$J_\zeta= \{b \in  \BB'_\zeta\mcolon b/I_\zeta \meet_{\BB''_\zeta}
  f_\epsilon(\gamma) = 0 \ \textrm{for every}\ \gamma < \lambda\}.$$ Then
   $J_\zeta$ is an ideal of $\BB'_\zeta$ extending $I_\zeta$ so $b
    \mapsto b/J_\zeta$ is a homomorphism. Further, $f_\epsilon(\gamma)$ is an
     atom of
 $\BB'_\zeta/J_\zeta$ for $\gamma< \lambda$.  These atoms are distinct and dense in $\BB'_\zeta/J_\zeta$. That is, $\BB_\epsilon$ is an atomic Boolean algebra.

\begin{snotation}\label{Beps}
Let $\BB_\epsilon$ be $\BB'_\zeta/J_\zeta$ with quotient map, $j_\epsilon(b) = b/J_\zeta$.
\end{snotation}

Now we  define $M_\epsilon$ by setting $P^{M_\epsilon}_1= \BB_\epsilon$
which contains $P^{M_\zeta}_1$;
$P^{M_\epsilon}_{4,1}$ is the injective image in $P^{M_\epsilon}_1$ of $P^{M_\zeta}_{4,1} \cup B_\epsilon$.
For $a \in P^{M_\epsilon}_{4,1}$ and $b \in P^{M_\epsilon}_{1}$,  set $R^{M_{\epsilon}}(a,b)$ if for some $\gamma$, $a = f_{\epsilon}(\gamma)/J_\zeta$ and
$f_{\epsilon}(\gamma)/J_\zeta \leq_{\BB_\epsilon} b /J_\zeta$.
Finally, let $D_\epsilon$ be the ultrafilter on $P^{M_\epsilon}_1$ generated by $$
D_\zeta \cup \{j_{\epsilon}(-b_{\zeta,\gamma})\mcolon \gamma < \lambda\} \cup \{j_{\epsilon}(-f_\epsilon(\gamma))\mcolon \gamma < \lambda\}.$$

We verify $M_\epsilon \in \bK_{-1}$ below. By Claim~\ref{succ}, we have the cardinality
 and atomicity conditions of
Specification~\ref{oplus}.(a) and (b); the definition of $I_\zeta$ guarantees,
 (c) and (d).(ii), (d).(iii).   The  elements $b_{\zeta,\gamma}$
 along with (our later) definition of $F^{M_\epsilon}_n(c)$ show d.i), $D_\epsilon
  \in \uf(M_\epsilon)$, (as no new $F_n(c)$ is in $ D_{\epsilon}$); the elements of $B_\epsilon$
   show $D_{\epsilon}$ is non-principal as each complement of an atom is in the ultrafilter.
 Note that  Specification~\ref{oplus}.(e) does not apply except in the $\omega+1$st stage of the construction.

%

For Specification~\ref{oplus} (f) (i), let $X\in \BB$ be a set of atoms of $M_\epsilon$
and note that we can choose
$b_X$ by conditions ii) and iii) in Definition~\ref{defI} of $I_\zeta$.


We can choose $P^{M_\epsilon}_2$ and
$F^{M_\epsilon}_n$ to satisfy Specification~\ref{oplus} (f) (ii). Fix an $\overline{A} \in \Ascr$ (as given by $\boxplus$).
 Fix a $c = c_
{\overline{A}}$ and define, using the last paragraph, the $F^{M_\epsilon}_n(c)$ as $b_{A_n}$, so that for each $n$,  $A_n
=\{\alpha<\lambda\mcolon f_\epsilon(\alpha) \leq_{P^{M_\epsilon}_1} F^{M_\epsilon}_n(c)\}$.  These are the only new $c\in P^{M_\epsilon}_2$.



Thus, it remains only to  show that
$M_\epsilon \in \bK_{-1}$. Most of the cases are obvious. E.g. for Definition~\ref{f1}.(8), just look at where the generators can be and recall countable free algebras are atomless.
Showing $M_\epsilon$ satisfies countable incompleteness, Definition~\ref{f1}.(7), is a bit more complex but we do so now.

($\blacklozenge$) If $a \in P^{M_\epsilon}_{4,1}$ and $c\in P^{M_\epsilon}_2$
 then $(\forall^{\infty} n)\ a\nleqslant_{M_\epsilon} F^{M_\epsilon}_n(c)$.

\noindent
If $c\in P^{M_\zeta}_2$, $F^{M_\epsilon}_n(c) = F^{M_\zeta}_n(c) \in P^{M_\zeta}_1$
 and we know by induction that $\blacklozenge$ holds for $a \in P^{M_\zeta}_{4,1}$.
 For $a\in P^{M_\epsilon}_{4,1}
 - P^{M_\zeta}_{4,1}$,
Definition~\ref{f1}.5, and
condition (vi   ) on $I_\zeta$  (from Definition~\ref{defI}) imply $a \leq_{M_\epsilon} b$ for every $b \in D_\zeta$.
 As $c \in P^{M_\zeta}_2$ and $D_\zeta \in \uf(M_\zeta)$, all but finitely many $n$,
 $e_n =F^{\zeta}_n(c)$, are {\em not} in $D_\zeta$.  So for all but finitely many $n$,  the complement
 $e_n^- \in D_\zeta$. That is, $a
  \leq_{M_\epsilon} e_n^-$; so  $a \wedge_{M_\epsilon} e_n = \emptyset$ as
   required.

If $c \in P^{M_\epsilon}_2-P^{M_\zeta}_2$ then
by our choice of $P^{M_\epsilon}_2$ and the $F^{M_\epsilon}_n$, there is an
$\overline{A}_c$ that is enumerated by the $F^{M_\epsilon}_n(c)$ and satisfies
$\blacklozenge$ by (i) of $\boxplus$ (Definition~\ref{boxplus}.(i)).
This completes the verification of $\blacklozenge$ at stage $\epsilon$
and so $M_\epsilon$ satisfies all the specifications of the induction.

%
%
%
%
%
%


%
%



%
%
%
%

\underline{case 4:}
$\epsilon = \omega+1$:

Only  clauses (c) and (e) of  Specification~\ref{oplus} are relevant.
Define  $P^{M_\epsilon}_2$ and $F^{M_\epsilon}_n$ to satisfy clause (e).
Since $P^{M_\epsilon}_i = P^{M_\omega}_i$ for $i = 0,1$, specification c) is immediate. 
This completes the construction.

\bigskip

{\bf The construction suffices.}

Having completed the induction, let $M= M_{\omega+1}$. Using specifications d) and a) of \ref{oplus}, it is straightforward to verify
that $M \in \bK_{-1}$ and the Boolean algebra is atomic.
 By \ref{clb}, $P^{M_\omega}_i$ for $i =0,1$ have cardinality $\lambda$.  And by \ref{clf},
the same holds for $P^{M_{\omega+1}}_2$.

We now show $M$ is nicely free.  Let $\bbar =\langle b'_\beta\colon \beta< \lambda 
\rangle$ enumerate $\langle  b_{n, \alpha}\colon n< \omega, \alpha < \lambda \rangle$
without repetition and such that $\{b_{-1,\alpha}:\alpha < \lambda\} = \{b'_{2\alpha  }: \alpha< \lambda\}$.
So this picks out a first level of $P^M_0$ which is enumerated by the $F^{M_0}_n(c)$
 for $c \in P^{M_0}_2$ and $n< \omega$ by case 1 of the construction.

Now, $\bbar$ satisfies the requirements in
Definition~\ref{b9} of nicely free. As, by Specifications~\ref{oplus}. \ref{clc},
\ref{cld}  and since $P^{M}_1$ is constructed as  the union of the
 $P^{M_n}_1$, $P^{M}_1/P^{M}_4$ is generated freely by $\bbar /P^{M}_4$.
Finally, clause c) of Definition~\ref{b9} holds by clause \ref{cld}.ii) of Specification~\ref{oplus}.



 The crux is to show $M= M_{\omega+1}$ is $P_0$-maximal. For this, assume for a contradiction:

(*)  $P^M_0$ is not maximal; by Lemma~\ref{f8}, there is a $D \in\uf(M_{\omega+1}) = \uf(M_{\omega})$.

\smallskip
\noindent
For every $n< \omega$,  is there a $d\in D$ such that $R(M_\omega,d) \cap M_n = \emptyset$?

\smallskip

$\underline{{\rm Possibility\ 1}:}$ For every $n< \omega$, the answer is yes, exemplified
 by $d_n\in D$. Now for each $a \in P^{M_n}_0$, $a \not \in R(M_\omega,d_m)$ for
  all $m \geq n$. So the sequence $\overline{d} = \langle d_n: n<\omega\rangle$
   satisfies the hypothesis of
Specification~\ref{oplus}.\ref{cle}  and so there is a $c \in P^{M}_2$ such that for
 each $n< \omega$, $F^M_n(c) = d_n$.  Thus, recalling Definition~\ref{f5}, $D \not \in \uf(M)$.

\smallskip
$\underline{{\rm Possibility\ 2}:}$ For some $n<\omega$, there is no such $d_n$;
 without loss of generality, assume $n>0$. We apply specification f) with $\epsilon = n$.
  Recall that $f_n$ is a $1$-$1$ map from $\lambda$ onto $P^{M_n}_{4,1}$. Let $g_1$
  be the following homomorphism from the Boolean algebra $P^{M_{\omega+1}}_1= P^{M_{\omega}}_1$
   into $\Pscr(\lambda)$: $g_1(b) = \{\alpha<\lambda: f_n(\alpha) \leq_{\BB_{M_\omega}}b\}$.  By Specification f.i) of
\ref{oplus}, the Boolean algebra $\BB$ provided by $\boxplus$ is contained in
the range of $g_1$.


Let $\Iscr_n$ denote the ideal of $P^M_1$ generated by $P^M_{4,1}- P^{M_n}_{4,1}$.
Since $D$ is non-principal, $\Iscr_n \cap D = \emptyset$. Now, $g_1$ maps any $b \in P^{M_\omega}_1
-P^{M_{\omega}}_4$ (and, thus, any $b\in P^{M_{\omega} }_1 -\Iscr_n$)
 to a nonempty subset of $\lambda$. Recalling $\Iscr_n \cap D = \emptyset$,
$D_1 = g_1(D)$ is an ultrafilter of the Boolean Algebra $\rg (g_1)$
and so $D_2 =D_1 \cap \BB$ is an ultrafilter of the Boolean algebra $\BB$. 
We show, $D_2$ is non-principal, i.e., for any $\alpha<\lambda$, $\{\alpha\} \not \in D_2$. As,
$f_n(\alpha) \in P^{M_\omega}_{4,1}$ and so $f_n(\alpha)$ is not in $D$.
 So $\{\alpha\} \not \in D_1$. Thus, $\lambda -\{\alpha\} \in D_1$ and
 so $\lambda -\{\alpha\} \in D_2$. So $\{\alpha\} \not \in D_2$ as promised.

Now we apply the second clause of $\boxplus$ to the ultrafilter $D_2$.
 Since we satisfied specification f.ii) in the construction, we can conclude there
 is $\overline{A} =\langle A_n \mcolon n< \omega\rangle \in \Ascr$ such that for
 infinitely many $k$, $A_k$ is in $D_2$.
Thus,  $u =\{k\colon A_k \in D\}$ is infinite.  We will finish the proof by
showing there is a $c$ such that $u=u_c$ (Definition~\ref{b9}) is the set of images of the $F^M_n(c)$.

Since we are in possibility 2
), each $A_k \in \BB$, $A_k \in \rg(g_1)$. So we can choose $d_k \in
P^{M_{\omega}}_1$ with $g_1(d _k) = A_k$. As $A_k \in D_2$, by the choice of
$D_1,D_2$ we have $d_k$ is in the ultrafilter $D$ from the hypothesis for contradiction: (*).

We show the sequence $\overline{d}= \langle d_k\mcolon k< \omega\rangle$
satisfies the hypothesis of clause e of Specification~\ref{oplus}.
First, $d_k \in P^{M_{\omega}}_1 - P^{M_{\omega}}_{4}$ as $D$ is a
 non-principal ultrafilter on $P^{M_{\omega}}_1$ so the first hypothesis is satisfied. Further, for every $a \in
 P^{M_{\omega}}_{0}$ all but finitely many $k$, $G^{M_{\omega}}_1 (a) \nleq_{M_{\omega}} d_k$ because $\overline{A} \in \Ascr$, which implies by
 $\boxplus$ ii) that for every $\alpha <
 \lambda$, for some $k_\alpha$, we have $k\geq k_\alpha$ implies $\alpha \not
 \in A_k$. Now by the definition of $g_1$, recalling $g_1(d_k) = A_k$, we
 have
$k\geq k_\alpha$ implies $f_k(\alpha) \nleqslant d_k$ (in $P^{M_{\omega}}_1)$.
So by Specification \ref{oplus}. f.ii), there is a $c \in P^{M_n}_2$ such that if
for all $k<\omega$, $F^{M_n}_k(c) = d_k$.  So, for each finite $k$,
 $d_k \in D$ and  $F^{M_{\omega+1}}_k(c) = d_k$. This contradicts $D \in
 \uf(M_{\omega+1})$ and we finish.   $\qed_{\ref{f11a}}$

%
%

\subsection{$\bK_1$ and $\bK_2$}\label{K1def}

We now introduce some terminology
 from \cite{BaldwinShelahhanfmax}. We first describe three subclasses of $\bK_{-1}$: $K^1_{<\aleph_0}$, the finitely generated models, their
 direct limits $\bK_1$
  and then the extension to $\bK_2$, the models of
 the complete sentence.

\begin{definition}[$\bK_{<\aleph_0}^1$ Defined]\label{k0}
  $M$ is in the class of structures $\bK^1_{<\aleph_0}$ if  $M \in \bK_{<\aleph_0}^{-1}$ and there is a {\em witness}
  $\langle n_*, \bB, b_*\rangle$
  such that:

  \begin{enumerate}
  \item $b_* \in P^M_1$ is the supremum of the finite joins of
      atoms in $P^M_1$. Further, for some $k$, $\bigcup_{j\leq k}
      P^M_{4,j} = \{c: c \leq b_*\}$ and for all $n>k$,
      $P^M_{4,n} =\emptyset$.
  \item $\bB = \langle B_n: n \geq n_*\rangle$ is an increasing
      sequence of finite Boolean subalgebras of $P^M_1$.


\item   $B_{n_*}\supsetneq    \{a\in P^M_1: a \leq b_*\} =
    P^M_4$;  the subset $$P^M_4
    \cup\{F^{M}_n(c)\colon n< n_*, c \in P^{M}_2\} $$  generates $B_{n_*}$.


Moreover, the Boolean algebra $B_{n_*}$ is free over the ideal
$P^M_4$  (equivalently, $B_{n_*}/P^M_4$ is a free Boolean
algebra\footnote{A further equivalence:
 $|Atom(B_{n_*})|/
|P^M_{4,1}|$ is a power of two.}).


  \item $\bigcup_{n\geq n_*}B_n = P^M_1$.
  \item $P^M_2$ is finite and not empty. Further, for each $c \in
      P^M_2$ the $F^M_n(c)$ for $n< \omega$ are independent over $P^M_4$.





\item  The set $\{F^M_m(c):m\geq n_*, c \in P^M_2\}$ (the
    enumeration is without repetition) is free from $B_{n_*}$
    over  $P^M_4$, 
     $B_{n_*} \supsetneq P^M_4$ and $F^M_m(c)
    \wedge b_* = 0$ for $m\geq n_*$. (In this definition, $0 =
    0^{P^M_1}$.)

 In detail,
let $\sigma(\ldots x_{c_i} \ldots)$ be a Boolean algebra term in
the variables $x_{c_i}$ (where the $c_i$ are in $P^M_2$  which is
not identically $0$. Then, for finitely many $n_i\geq n_*$ and a finite sequence of  $c_i\in P^M_2$:
$$ \sigma(\ldots F^M_{n_i}(c_i) \ldots) > 0.$$
Further, for any non-zero $d \in B_{n_*}$ 
with $d \wedge b_* =0$,  (i.e. $d \in B_n - P^4_M$),
$$
\sigma(\ldots F^M_{n_i}(c_i) \ldots) \wedge d > 0.$$

\item For every $n\geq n_*$, $B_n$ is generated by $B_{n_*} \cup \{F^M_m(c):n > m\geq n_*, c \in
P^M_2\}$. Thus $P^M_1$ and so $M$ is generated by $B_{n_*}\cup
P^M_2$.

\end{enumerate}

\end{definition}

Recall some terminology from \cite{BaldwinShelahhanfmax}.

\begin{definition}[$\bK_1,\bK_2$ Defined]\label{richname}\begin{enumerate}
\item $\bK_1$ denotes the collection of all direct limits of models in $\bK^1_{<\aleph_0}$.
     \item We say a model $M$ in $\bK_1$ is rich  if for any $N_1, N_2 \in \bK^1_{<\aleph_0}$ with $N_1 \subseteq N_2$ and $N_1 \subseteq M$, there
is an embedding of $N_2$ into $M$ over $N_1$.
\item $\bK_2\subseteq \bK_1$ is the class of rich models.
\end{enumerate}
\end{definition}

Note that the free generation in item 6 of Definition~\ref{k0} is not preserved
by arbitrary direct limits and so is not a property of each model in $\bK_1$.
In particular, as $M_*$ is corrected to a model of $\bK_1$, we check the freeness
only for finitely generated submodels as it will be false in general.

Since $\bK^1_{<\aleph_0}$ has joint embedding, amalgamation and only
countably many finitely generated models, we construct in the usual
way a generic model; thus $\bK_2$ is not empty.

\begin{fact}\label{getgen}  There is a countable generic model $M$
for $\bK_1$ (Corollary~3.2.18 of \cite{BaldwinShelahhanfmax}). We denote its
Scott sentence by $\phi$. $\bK_2$ is the class of models of this $\phi$.
\end{fact}

\section{Correcting $M_*$ to a model of $\bK_2$}\label{corr}

\numberwithin{theorem}{section}

We now `correct' the $P_0$-maximal model of $\bK_{-1}$, $M_*$, constructed in Section~\ref{1approx},
  to obtain a $P_0$-maximal model $M$ (Definition~\ref{f2*}) of the complete sentence
   constructed in \cite{BaldwinShelahhanfmax}, i.e. $M \in \bK_2$. In Theorem~\ref{f56} we
 modify $M_*$,
 to construct a model $M\in \bK_2$ with $P^M_2 \subseteq P^{M_*}$ and
 redefining the $F_n$, but
   retaining  $M \myrestriction (P_0^M \cup P_1^M) = M_* \myrestriction (P_0^{M_*} \cup P_1^{M_*})$.
     The old values of $F^{M_*}_n$
will be used to divide the work of ensuring each ultrafilter $D$ is not
in $\uf(M)$ by for each $D$,  attending one by one to only those
 $c$ with infinitely many $F^{M_*}_n(c)$ in $D$.

We now  describe some of the salient properties of the model $M$ obtained by `correcting' the $M_*$ of Section~\ref{1approx}.

\begin{remark}[The Corrections]\label{correnum2}

\begin{enumerate}
\item The domains of the structures constructed in this section are subsets of $M_*$;  the $F_n$ are redefined so the new structures are substructures   only of the reduct of $M_*$ to $\tau -\{F_n: n<\omega\}$.

    \item In particular, for all the $M$ considered in Section~\ref{corr}, $P^M_1 =P^{M_*}_1$
     and these Boolean algebras have the same set of ultrafilters.
       However, $\uf(M) \neq \uf(M_*)$ as the  definition  of $\uf$ depends on
        properties of the  $F_n$.

        \item The set $\{F^M_n(c): c\in P^{M}_2\}$ is not required to be an independent subset to
        put  $M \in \bK_{-1}$.
    \item Lemma~\ref{f53} demands a sequence of finite Boolean algebras $B_n$ to
     witness finitely generated substructures belong to $\bK_1$ (not required for $\bK_{-1}$).  The stronger class of $\bK_1$-free structures \cite[Definition 3.2.11]{BaldwinShelahhanfmax}, which is closed under extension by members of $\bK_1$ and so has no maximal models plays no active role in this paper.
     In particular, the final counterexample, Theorem~\ref{f56}, is in $\bK_1$ but is not $\bK_1$-free.
\item The proof is in ZFC. The proof in \cite{BaldwinShelahhanfmax} that a non-maximal  model in $\lambda$ makes $\lambda$ measurable depends on $\diamond$.
    \end{enumerate}
    \end{remark}

The main task of this section is to prove:
 \begin{theorem}\label{realthm}
 If $\lambda$ is less than the first measurable cardinal,  $2^{\aleph_0}< \lambda$, and
  for some $\mu$, $2^{\mu} =\lambda$ (whence $\lambda^\omega =\lambda$),
  then there is a $P_0$-maximal model in $\bK_2$ of cardinality $\lambda$.
 \end{theorem}

Conclusion~\ref{hyp},
 summarises the results of the construction in 
Theorem~\ref{f11a}, specifically to fix our assumptions for this section.

\begin{conclusion}\label{hyp} If
$\lambda$ is as in Theorem~\ref{realthm} then there is
 a model $M_*$ with  $|M_*|= \lambda$ satisfying:
\begin{enumerate}
    \item $P^{M_*}_1$ is an atomic Boolean algebra and $M_*$ is $P_0$-maximal.  Further,
    $|P^{M_*}_i| = \lambda$ for $i =0,1$.
    \item $P^{M_*}_{4,1}$ is the set of atoms of $M_*$.
    \item $M_*$ is nicely free (Definition~\ref{b9}); in particular, $P^{M_*}_{1}/P^{M_*}_{4}$
    is a free Boolean algebra of cardinality $\lambda$.
 \end{enumerate}
 \end{conclusion}

 In order to `correct' $M_*$ to a model in $\bK_2$, we lay out some notation for
  the indexing of the tasks performed in the construction, the generating set of $P^{M_*}_1$,  and the free basis of the
 Boolean algebra $P^{M_*}_1/P^{M_*}_4$.

  \begin{notation}\label{f33} We define a family  of trees of sequences:
  \begin{enumerate}
  \item
For $\alpha < \lambda$, let    $\Tscr_\alpha= \{\langle \rangle\} \cup \{\alpha \add \eta; \eta \in {}^{<\omega}3\} $  and $\Tscr  = \bigcup_{\alpha<\lambda}\Tscr_\alpha$.
 \item  $\lim (\Tscr_\alpha)$ is the collection of paths through $\Tscr_\alpha$.
     \end{enumerate}
 \end{notation}

 Combining the requirements for constructing $M_*$ (Specification~\ref{oplus})
 and the Definition~\ref{b9} of nicely free, we have

\begin{claim}[Fixing Notation]  \label{f34} Since $M_*$ is nicely free,
 without loss of generality, we may assume:

\begin{enumerate}
\item The universe of  $M_*$ is $\lambda$ and the $0$ of $P^{M_*}_1$ is the ordinal $0$.
    \item We can choose
sequences of  elements of $P^{M_*}_1$, $\bbar =\langle b_\eta \colon \eta \in \Tscr\rangle$
 so that their images in the natural projection of $P_1^{M_*}$ on $P_1^{M_*}/P_4^{M_*}$  freely  generate
  $P_1^{M_*}/P_4^{M_*}$.

  \item For every $a\in P^{M_*}_{4,1}$ and the even ordinals $\alpha<\lambda$,  there is an $n$
 such that   for any $\nu \in \Tscr_\alpha$,
  $\lg(\nu) \geq n$ implies
  $a \wedge b_\nu  = 0$.
\end{enumerate}
\end{claim}


Proof.  The only difficulty is deducing   from c) of Definition~\ref{b9}
(nicely free) that 3) holds.  For that, we can insist
 that for each even $\alpha$, for some $c\in P^{M_*}_2$,
$\{b'_{\omega  \alpha + n}\colon n< \omega\}$ enumerates
 $u_c = \{F^{M_*}_n(c): n< \omega\}$ (from Definition~\ref{b9}.c).
Now for $\alpha > 0$, let $\langle b_\eta \colon \eta \in
\Tscr_\alpha \setminus \{\langle \rangle\}\rangle$
list $\{b'_{\omega  \alpha + n}\colon n< \omega\}$ without repetition
and $\langle b_\eta \colon \eta \in \Tscr_0 \rangle $ list  $\{b'_n\colon n< \omega\}$.
By Definition~\ref{f1}.7 ($\bK_{-1}$) we have:
for every $a\in P^{M_*}_{4,1}$ for all but finitely many $n$,
$a \wedge b'_{\omega  \alpha + n} = 0_{P^{M_*}_{1}}$;
 whence  for even $\alpha$ all but finitely many of the
 $\nu \in \Tscr_\alpha$  
 satisfy $a \wedge b_\nu = 0_{P^{M_*}_{1}}$.
 $\qed_{\ref{f34}}$




\medskip
Note that Claim~\ref{f34} provides a  $1$-$1$ map from $P^{M_*}_2$ to ordinals less than
$\lambda$. We introduce the collection of models that is the starting point for the following construction.

  \begin{definition}[$\MM_1$ Defined]\label{f37} Let
   $\Mscr_1 = \Mscr_1(\lambda)$ be the set of $M\in \bK_{-1}$ such that the
    universe of $M$ is contained in $\lambda$, the universe of $M_*$, and for
     $i< 2$, (or $i =4$ or $(4,1)$) $P_i^M = P_i^{M_*}$,  $M \myrestriction (P_0^M \cup P_1^M) =
      M_* \myrestriction (P_0^{M_*} \cup P_1^{M_*}$) while $P^M_2 $ will not equal $P^{M_*}_2$.

  \end{definition}

  The posited $M_*$ differs from any $M\in \Mscr_1$ only in that $P_2^M$
  is, in general, a proper subset of $P_2^{M*}$ and the newly defined
  $F_n^M(c)$ (usually) do not equal the $F_n^{M_*}(c)$.
    We now spell out the tasks
  which must be completed to correct $M_*$ to the required member of $\bK_2$.
  The $F^{M_*}_n(c)$ are used  as   oracles.

   \begin{definition}[Tasks]\label{f39} \begin{enumerate}
   \item Let $\bT_1$, the set of 1-tasks, be the set of pairs $(N_1,N_2)$ such that:
      \begin{enumerate}
   \item $N_1 \subseteq N_2 \subseteq \lambda$

   \item $N_1, N_2 \in \bK^{1}_{<\aleph_0}$

   \item $N_1 \subset M$ for some $M \in \Mscr_1$. More explicitly, $P^M_2 \subseteq P^{M_*}_2$
   and\\ $N_1\myrestriction (P^M_0 \cup P^M_1) \subseteq M_*$ and
   $(F^M_n \myrestriction P^{N_1}_2 ) = F^{N_1}_n$ for each $n$.
   \end{enumerate}

   \item Let $\bT_2$, the set of 2-tasks, be the set of $c \in P_2^{M_*}$.
   \item $\bT= \bT_1 \cup \bT_2$.
   \item Let $\langle \tbar_\alpha\colon \alpha < \lambda\rangle$ enumerate $\bT$.
   \end{enumerate}
   \end{definition}

   Note $|\bT_1| = |\bT_2| =|\bT|$.

   \begin{definition}[Task Satisfaction]\label{f41}
   The task $\tbar$ is {\em relevant} to the structure $M$ if
   $M\in \MM_1$ and i) if $\tbar$ is 1-task $(N_1,N_2)$ and $N_1 \subseteq M$ or ii) if $\tbar$ is a 2-task $c$ and $c \in P_2^{M}$.

   We say $M \in \Mscr_1$ {\em satisfies the task} $\tbar$ if either:
   \begin{enumerate}[A)]
   \item $\tbar = (N_1,N_2)\in \bT_1$ (so $N_1 \subset M$) and there exists
   an embedding of $N_2$ into $M$ over $N_1$.
   \item $\tbar =  c$, where $c \in P_2^{M_*}$, is in $\bT_2$ and for every ultrafilter $D$ on $P_1^{M}$,
        such that for infinitely many $n$,
       $F_n^{M_*}(c) \in D$,
       there is a $d\in   P_2^{M}$ such that for infinitely many $n$,
       $F_n^{M}(d) \in D$.
%
%
%
%

       \end{enumerate}
       \end{definition}

Recall Definition~\ref{f5} of $\uf(M)$ and Lemma~\ref{f8} connecting $\uf(M)$  with $P_0$-maximality of $M$.


       \begin{claim}\label{f44}

           If $M \in \Mscr_1$ satisfies all tasks in $\bT$ and is in $\bK_1$ then from
           satisfying the $\bT_2$ tasks, $M$ is $P_0$-maximal   and satisfying the tasks in $\bT_1$ guarantees it is in $\bK_2$.

               \end{claim}

               Proof.
            For  $P_0$-maximality of $M$, it suffices, by Lemma~\ref{f8} (since $\MM_1 \subseteq \bK_{-1}$), to show $\uf(M) = \emptyset$.
              But, since $\uf(M_*) =\emptyset$, for every ultrafilter $D$ on $P^{M_*}_1$
             there is $c\in P^{M_*}_2$  with $S^{M*}_c(D)$
             infinite (Definition~\ref{f5}); satisfying task $c$ means there is $d\in P^{M}_2$
             such that $S^{M}_d(D)$ is infinite and so $D$ is not in $\uf(M)$.
             Since $M$ and $M^*$ have the same ultrafilters, this implies
             $\uf(M) = \emptyset$, as required.  Since we have assumed $M\in \bK_1$, the second
             assertion follows by realizing that satisfying all the tasks
             in $\bT_1$ establishes the model is rich, which suffices by Fact~\ref{getgen}. $\qed_{\ref{f44}}$.


\bigskip

   Definition~\ref{f50} lays out the use of the generating elements $b_\eta$ in correcting the $F^{M*}_n$ to require independence while maintaining that
   infinite intersections of members of the ultrafilter under consideration are empty.   The
   infinite sequence $\eta_d$ will guide the choice of
    $F^M_n(d)$.

%


The following facts about the relation of symmetric difference and ultrafilters are central for calculations below.
\begin{remark}\label{backgrbauf} Recall that the operation of symmetric difference is associative.

  \begin{enumerate}
  \item 
  Suppose $\BB_1 \subseteq \BB_2$ are Boolean algebras with $a \in \BB_1$, and
  $b_1 \neq c_1$ are in $\BB_2$, and $\{b_1,c_1\}$ is
  independent over $\BB_1$ in $\BB_2$.

   The element $(b_1\vartriangle c_1) \vartriangle a\in \BB_2$ is independent over $\BB_1$.
   More generally, if $\{b_i,c_i:i<\omega\}$ are  independent over $\BB_1$,
   $\{a_i:i< \omega\} \subseteq \BB_1$ ,
   $e_i = b_i  \vartriangle c_i \vartriangle a_i   $,e and $f_i = b_i \vartriangle c_i$
    then each of $\{e_i: i< \omega\}$ and $\{f_i: i< \omega\}$ are independent over $\BB_1$.


\item 
%
%
  Let  $D$ be an ultrafilter on a Boolean algebra $\BB$.

  \begin{enumerate}
  \item
 For $a,b \in D$,

   ($a \in D$ iff $b \in D$) if and only if $a\vartriangle b\not \in D$.

 \item                 If $a_0, a_1, a_2 \in \BB$ are distinct  then   at least one of $a_i \vartriangle a_j \not \in D$.

                 \item  More importantly for our use later, it is easy to check: \

                        ($a_0 \in D$  iff $a_1 \in D$)

                         iff $$(a_0 \vartriangle a_1 \vartriangle a_2) \in D \leftrightarrow a_2 \in D.$$

 \end{enumerate}
  \item If $a$ is an atom, $a \wedge b_0 =0$   and $a \wedge b_1 = 0$,
  then $a \wedge (b_0 \vartriangle b_1) = 0$.


\end{enumerate}
\end{remark}

Proof. 1)   If  the element $(b\vartriangle c) \vartriangle a \in \BB_2$ is not independent over $\BB_1$ there is a polynomial $p$ over $\BB_1$ with
$p((b\vartriangle c) \vartriangle a)\in \BB_1$. But then, by
 Observation~\ref{obs1a}, $p(x,y) = p((x\vartriangle y) \vartriangle a)$
  is also a polynomial over $\BB_1$ witnessing $\{b,c\}$ is dependent over $\BB_1$.
In the more general case any polynomial witnessing dependence in $n$ of the $e_i$ ($f_i$) give a
polynomial in $2n$ of the $a_i, b_i,c_i$ witnessing dependence of the original set.

2) For a), if, say $a\in D$ and $b\not \in D$, then $a-b$ and hence
$a \vartriangle b \in D$ so we have `left to right' by contraposition.
 If both are in $D$,  so is their meet which is disjoint from $a \vartriangle b$ so
  $a \vartriangle b \not \in D$.  Since $a^- \vartriangle b^- = a \vartriangle b $,
   we have the result if neither is in $D$.

b) holds since  the intersection over all pairs $i,j < 3$ of the $a_i \vartriangle a_j$ is empty.
  And c) is propositional logic from a) and b).

3) $a \leq (b_0^{-} \wedge  b_1^{-}) \leq (b_0^{-} \vartriangle  b_1^{-}) \leq (b_0  \vartriangle  b_1)^{-}$. As $a$ is an atom, $a \wedge (b_0  \vartriangle  b_1) =0$.
$\qed_{\ref{backgrbauf}}$

\medskip

    We define a class   $\Mscr_2\subseteq \MM_1$ such that for each $d\in P^M_2 \in \Mscr_2$ there is an ordinal $\alpha_d$, a
     tree of elements of
     $P^M_1$, indexed by sequences in  $(\Tscr_{\alpha_d})\subseteq
     {}^{<\omega}3$, a target path $\eta_d$ through that tree and a sequence
      $a_{d,n}$, whose indices are not in $\Tscr_{\alpha_d}$, but which
      satisfy that each $a \in P^{M_*}_{4,1} =P^{M}_{4,1}$ is in at most
      finitely many $a_{d,n}$.
    In the construction (Theorem~\ref{f56}) of a model in $\MM_2$, $\eta_d$  guides
    definition of the
    sequence $F^M_k(b_{\eta_d})$.  The $a_{d,n}$ are introduced
     to make Definition~\ref{f50}.B
    uniform.  In cases 2 and 3 of Theorem~\ref{f56} $a_{d,n}$ is always $0$.
    In case 4, where the $F^M_n(d)$
    are defined as $M$ is corrected from $M_*$, $a_{d,n} = F^{M_*}_n(d)$.
       The result is the values of the $F^M_{n}(d)$
     are both independent over a finite initial segment and satisfy
      $\bigwedge_{n<\omega} F^M_n(d) = \emptyset$.
The next definition abstracts from this construction to identify the key ideas of the
proof that if $M\in\MM_2$ then $M \in \bK_1$ (Lemma~\ref{f53}) and further that
there are $M \in \MM_2$ that are in $\bK_2$.  The notation
$\langle Z \rangle$ denotes the Boolean subalgebra of $P^M_1$ generated by $Z$.

 \begin{definition}[$\MM_2$ Defined]\label{f50} Let $\Mscr_2$ be the set of
 $M \in \MM_1$  such that there is a sequence $\wbar =\langle (\alpha_d,\eta_d,
  a_{d,n})\colon d \in P_2^M, n< \omega\rangle$ witnessing the membership,
  which means:

              \begin{enumerate}[A]
         \item {}      \begin{enumerate}
     \item   For each $d\in P^M_2$,      $\alpha_d < \lambda$ is even and $d_1 \neq d_2$ implies
$\eta_{d_1}\neq \eta_{d_2}$.  (In case 4 of Lemma~\ref{f56}, many $d_\eta$ have the same $\alpha_\eta$.)

\item $\langle \alpha_{d}\rangle \lhd\eta_{d} \in \lim(\Tscr_{\alpha_d})$.
\end{enumerate}

  \item  For each $n< \omega$, there are\footnote{In applications, the $a_{d,n}$  are either
  $0$ or $F^{M_*}_n(c)$ (for an appropriate $c \in P^{M_*}_2$).}  $a_{d,n}$  in
   $P^{M_*}_1 =P^M_1$  such that for
   each $d \in P^M_2$, there are 
   distinct\footnote{I.e., $\nu_1[d,n]$ depends on $d$
   and $n$.} $\nu_1[d,n]$
   and $\nu_2[d,n]$ that extend
    $\eta_d\myrestriction n$,  $\nu_i(0) =\alpha_d$, and have length $n+1$
   such that:
      \begin{enumerate}
     \item   For every $n$,
             $$F^M_n(d) = (b_{\nu_1[d,n]} \vartriangle b_{\nu_2[d,n]})\vartriangle  a_{d,n}.$$
%




             \item for each $a \in  P^{M_*}_{4,1}$ and each $ d \in P_2^M$,
             there are only finitely many $n$ with $a \leq_{P^{M_*}_1} a_{d,n}$.
         \end{enumerate}

         \item ($k^1_Y$) For each finite $Y \subseteq P^M_2$ there is a list
         $\langle d_\ell: \ell < |Y|\}$ of $Y$ such that:


              \begin{enumerate}
              \item The $d_\ell$ list $Y$ without repetition
              and
               $\alpha_\ell = \alpha_{d_\ell}$.
              \item   If  $i_1 < i_2 < i_3 < |Y|$ and $\alpha_{i_1} = \alpha_{i_3}$
               then $\alpha_{i_2} = \alpha_{i_1}$.

          \item     There is\footnote{See proof of {\em goal} in Lemma~\ref{f56}.} a $k_1 = k^Y_1$ such that

          \begin{enumerate} \item For $i\neq j$, both less than $n$,
     $\eta_i \myrestriction k^Y_1 \neq  \eta_j \myrestriction k^Y_1$.

      \item   
       Set
            $ W\subseteq P^{M_*}_1$ as:
 \begin{eqnarray}
 W & = & \{a_{d_k,n}: k<|Y| \wedge n
             <\omega\}\\
   \nonumber  & &          \cup \  \{F_i^M(d_k): k< |Y|, i< k^Y_1\}.
 \end{eqnarray}



    Then $ W$ is included in
 the  subalgebra 
 $\BB^0_Y$  of $P^{M}_1$ generated by
  $$\{b_\nu\colon \bigwedge_{i<|Y|} (\eta_i \myrestriction k^Y_1) \ntrianglelefteq
      \nu
      \} \cup \{b_{\langle \rangle}\} \cup P^M_{4,1}$$
        where $\eta_i$ abbreviates $\eta_{d_i}$.
\end{enumerate}



              \end{enumerate}
 \end{enumerate}
\end{definition}

              Note that the $B^0_Y$ is a cocountable subset  of $P^M_1$
              (the complement is contained in  finite set of countable trees).

We will apply the following lemma three times to show that for $M\in \MM_2$, the set $\{F^M_n(c)\}$
is countably incomplete (witnessing Definition~\ref{f1}.7). It is a straightforward application of
of Remark~\ref{backgrbauf} to Definition~\ref{f50}.2.

\begin{lemma}\label{countinc}
Let $M \in \MM_1$.  For any $\langle \alpha_d, \eta_d, a_{d,n}\rangle$ as in Definition~\ref{f50}, (in particular $\alpha_d$
is even) and any atom $a\in P^{M_*}_{4,1}$,
%
%
for all but finitely many $n$
$$a \wedge (b_\nu \vartriangle  b_\rho  \vartriangle a_{d,n}) = 0.$$
\end{lemma}
Proof. 
Recall from  \ref{f34}.3,  that for
every $a\in P^{M_*}_{4,1}$ and the
 even ordinals $\alpha<\lambda$,
 there is an
$n$,  such that for any $\nu, \rho \in \Tscr_\alpha$ with $\lg(\nu) \geq n$
 and $\lg(\rho) \geq n$,
 $a \wedge b_\nu =0$ and $a \wedge b_\rho =0$.
 Definition~\ref{f50}.B.b asserts each $d$ and for sufficiently large $n$, $a \wedge a_{d,n} =0$.
Apply Remark~\ref{backgrbauf}.3 twice.
 $\qed_{\ref{countinc}}$

We will show in Lemma~\ref{f53} that members of $\MM_2$ are in $\bK_1$ and
 then in Theorem~\ref{f56} that there are structures in $\MM_2$ that are in $\bK_2$.
Two main features distinguish $\bK_1$ from $\bK_{-1}$.
The $F_n(d)$ retain the `countable incompleteness'
property from $\bK_{-1}$ but
also must be independent; $M \in \bK_1$ when $M$ is a direct limit of members of $\bK^1_{<\aleph_0}$.

    \begin{lemma}\label{f53} If $M \in \MM_2$, then $M \in \bK_1$.
    \end{lemma}

    Proof. Suppose $M \in \MM_2$. Let $Y \subset P_2^M$ and $X \subset P_1^M$ be
    finite; 
    we shall find $N =N_{XY}\in \bK_{<\aleph_0}^1$ such
    that $Y \cup X \subseteq N \subseteq M$; this suffices.
As, 
$\bK_1$ is defined to be the collection of direct limits of
 finitely generated structures\footnote{The proof of Lemma{f53} shows there
  is a common substructure of $M$ containing  any finite collection
 of finitely generated (as in this argument) substructures of $M$.} in $\bK_{<\aleph_0}^1$.

 Our two main jobs in proving Lemma~\ref{f53} are to find an $N, n_*, b_*$  in which
 \begin{enumerate} \item
   the $F^M_k \myrestriction N$ satisfy
 property 6 (independence) of Definition~\ref{k0} over a $B_{n_*}$ and property 7 of Definition~\ref{f1} and then
  \item construct $N =\bigcup_{n< \omega} B_n$ for  finite
 Boolean algebras $\langle B_{n}:n\geq n_*\rangle$  that
  witness 2 and 3 of Definition~\ref{k0}.
  \end{enumerate}

 The finite $k_1=k^{Y}_1$ specified in Definition~\ref{f50} depends only
 on $Y$; in the next definition we increase $k_1$ to a $k^{X}_1 = k^{XY}_1$  and using
 the definition of $\MM_2$ show the $F^M_k(d)$ are independent over $X$ for $k \geq k^{XY}_1$.
We need $k^{XY}_1$ only to prove Lemma~\ref{f53}.

 We build two increasing chains of length $|Y|$ of subboolean algebras satisfying
 conditions described in Definition~\ref{k1req}.
The $\BB^\ell_{XY}$ will be cocountable, while the $\FF_\ell$ will be
  countable. The existence of $k^{XY}_1$ satisfying the conditions of
   Definition~\ref{k1req} is proved in Fact~\ref{basecontain}.


\begin{definition} [$k^1_{XY}$] \label{k1req}  Let the sequence
 $\langle (\alpha_d,\eta_d,  a_{d,k})\colon d \in P_2^M, k<\omega\rangle$
 {\em witness} $M\in \MM_2$ as in Definition~\ref{f50}.
    Let $X\subset  P^M_1$ and
 $\langle d_i: i<n\rangle$ enumerate $Y \subset  P^M_2$ without repetition and denote,
    for $i<n$, $\eta_{d_i}$ by $\eta_i$ and $\alpha_{d_i}$ by $\alpha_i$.
     Without loss, the  $\langle \eta_i(0)\colon  i<n\rangle$ are non-decreasing;
  \begin{enumerate}[A]
     \item
    Fix $k_1 = k^{XY}_1$ such that
     \begin{enumerate}
     \item $k^{XY}_1 \geq k^{Y}_1$ (see Definition~\ref{f50}.B);
      \item   $\langle \eta_i\myrestriction k^{XY}_1\colon i<n\rangle$ are distinct for  $i<n$;
\item $k^{XY}_1 \geq \max \{ \lg(\nu): b_\nu \in \langle X \cup \{F^M_k(d_i):i<|Y|\}\rangle, k< k^Y_1\}$.

\end{enumerate}
\item We consider the following sets determined by $X \cup Y$ and the $\eta_i$.
 \begin{enumerate}
     \item

$\bF_{\leq 0} =\bF_0 = X\cup \{F^M_k(d_i) \colon i< |Y|, k \leq  k^{XY}_1\}$;


\item For $1\leq \ell < |Y|$,
$\bF_\ell = \{F^M_k(d_\ell):k \geq k^{XY}_1\}$;
\item  $\bF_{\leq\ell+1} = \bF_{\leq\ell} \cup\bF_{\ell+1}$ ;
\item $\FF^\ell = \langle \bF_{\leq \ell}\rangle_{M}$.
%

\end{enumerate}
\item $$\BB^\ell_{XY} = \{b_\nu\colon \bigwedge_{\ell < i< n} (\eta_i \myrestriction k^{XY}_1) \ntrianglelefteq
     \nu \textrm{ for } i< \ell+1\} \cup \{b_{\langle \rangle}\} \cup P^M_{4,1}.$$

\end{enumerate}
\end{definition}

For each $\ell$, $\BB^\ell_{XY} \supseteq \BB^\ell_{Y}$ since $k^{XY}_1 \geq k^{Y}_1$
and $\BB^{\ell +1}_{XY} \supseteq \BB^\ell_{XY}$. In the proof of Lemma~\ref{f53.5}
$B_{n_*}$ will be $\FF^0$ and $N$ will be
$\FF^{n-1}$.

%

%
%

Since $X$ and $Y$ are finite we now choose $k^{XY}_1$ to satisfy conditions 1-3 of Definition~\ref{k1req};
we now show the other conditions are satisfied.

\begin{fact}\label{basecontain} There is a $k_1 = k^{XY}_1$ such
 that for each $\ell$, $\bF_\ell$ is contained in $\BB^\ell_{XY}$.
\end{fact}
Proof. Recall (Claim~\ref{f34}) that $M_*$ is free on the
$\{b_\eta: \eta \in \Tscr\}$ modulo the $P^{M_*}_4$.
Choose $k^{XY}_1$ larger than the length of any $\nu$ such that for some $x\in X$, $b_\nu$ is
a generator in a minimal representation of $x$ or $\nu(0)
 \in \overline{\alpha}=\{\alpha_0, \ldots \alpha_{n-1}\}$.
Then,

$$\bF_0 \subseteq \langle\{b_\nu: \nu \in \Tscr, \lg(\nu)< k^{XY}_1\}\rangle
 \cup \{b_{\langle\rangle}\} \cup P^M_4
\subseteq \BB^0_{XY}.$$
Recall from Definition~\ref{f50}.D, that as $\ell$ increases $F^M_k(d_i)$ for $i< \ell$ and all $k$ are admitted
to $\BB^\ell_{XY}$ and so $\bF_\ell \subseteq \BB^{\ell}_{XY}$.
  $\qed_{\ref{basecontain}}$




\smallskip
      To establish job 1) we need the following claim.

\begin{lemma}\label{f53.5} For each $1\leq \ell< n$, $\bF_\ell$ is independent over $\BB^0_{XY}$ mod $P^M_4$.
\end{lemma}

%
%
%
%
%
   Proof.
   We prove this claim by showing by induction on $\ell\leq |Y|=n$:
 $$(\oplus_\ell) \ \ \bF_\ell = \{  F_k^M(d_i)\colon k\geq k^{XY}_1 \mbox{ \rm and } i<\ell\}$$ is independent in $P^M_1$ over
 $\BB^{\ell-1}_{XY}$ mod $P^M_4$.

 For
 $1 \leq \ell <|Y|$, the induction on $\ell$ shows incrementally, at
 stage $\ell+1$, the independence of the $b_{\eta_\ell\myrestriction r}$ with $r \geq k^{XY}_1$
 over $\BB^\ell_{XY}$.
By Claim~\ref{f34}.2 and the choice of $r \geq k^{XY}_1$, the $\{b_{\nu_1[d_\ell,r]}:r \geq k^{XY}_1\}$
are independent mod $P^M_4$.
Thus  (using the $f_i$ from Remark~\ref{backgrbauf})
     the infinite set
     $ \{b_{\nu_1[d_\ell,n]} \vartriangle b_{\nu_2[d_\ell,n]}) \colon i\in \{0,1\}, n\geq k^{XY}_1\} $ is
    independent over $\BB^{\ell-1}_Y$.
By Definition~\ref{f50}.C) the  $\{a_{d_\ell,k}\colon k\geq k^Y_1\}$ are in $\BB^0_Y\subseteq \BB^\ell_{XY}$. Further, by
Definition~\ref{f50}.B) for all $n$: $$F^M_n(d_\ell)=
  (b_{\nu_1[d_\ell,n]} \vartriangle b_{\nu_2[d_\ell,n]})\vartriangle  a_{d_\ell,n}.$$
  So,
 Lemma~\ref{backgrbauf}.2 (now using the $e_i$)  implies $\bF_\ell$ is independent over  $\BB^\ell_Y$.
 Since independence is transitive (Lemma~\ref{obs1}.3)
    $\bF_\ell$ is independent over $\BB^0_Y$.
 $\qed_{\ref{f53.5}}$

We continue the proof of Lemma~\ref{f53}.
By Lemma~\ref{countinc},
   for sufficiently large $n$,  $a\nleqslant F_n^M(d_\ell)$.  So the countable incompleteness
   condition in the definition of $\bK_{-1}$ is satisfied.
This completes job 1). To accomplish job 2) and finish the proof of Lemma~\ref{f53}
 by satisfying conditions 2-4 of Definition~\ref{k0}, we must define appropriate $P^N_i$
 and find a sequence of finite Boolean
 algebras $B_n$ witnessing that $N \in \bK^1_{<\aleph_0}$.
%
Let $P^N_1 = \FF^{n-1}$.
We have  $P^N_1$ is freely generated (modulo the ideal generated by the atoms of $B_{n_*}$)
 by the countable set $\bF_{|Y|}$ over $B_{n_*}= \FF^0$.
  Let $b_*$ be the supremum of the atoms in $B_{n_*}$, and $P^N_{4}$ the predecessors of $b_*$.

 For $m\geq n_*$, let $B_m$ be generated by $B_{n_*}$ and the first $m$ elements
 of this generating set.   Then, $P^N_1 = \bigcup_{n_*\leq m <\omega}B_m $ and $P^N_1/P^N_4$ is atomless.
 Set\footnote{$G^M_1$ is from Definition~\ref{f1}.5.} $P^N_2 = Y$  and  $P^N_0 = \{(G^M_1)^{-1}(a): a\in P^M_{4,1} \cap P^N_1\}$;
 thus $P^N_{4,1} \subseteq B_{n_*}$.
 Boolean algebras are locally finite and we can recognize whether
  $\langle X\rangle$ is free if by whether it has $2^{|X|}$ atoms.
Thus, we can refine the sequence $B_m$ to finite free algebras to witness  that
$N \in \bK^1_{<\aleph_0}$. Since $X$ and $Y$ were arbitrary,
$M\in \bK_1$.
 $\qed_{\ref{f53}}$

    \medskip

   This completes the proof of Lemma~\ref{f53}. Now we show $\MM_2$ is non-empty and at least one member satisfies all the tasks.
     In case 4) of this argument we address the requirement that
      $\uf(M_\alpha) = \emptyset$ and so $\uf(M) = \emptyset$ as well. We need the following
      observation because as the construction proceeds, an $N_1$ may become a substructure of $M_\beta$
      because some value of an $F_n$ is newly defined on a point of $P^{M_\beta}_2$.

      \begin{notation}\label{priority} 
    We can enumerate $\bT$ as $\langle t_\alpha\colon \alpha< \lambda\rangle$ such
    that each task appears $\lambda$ times, as we assumed
    in Hypothesis~\ref{hyp} that $ \lambda=  \lambda^{\aleph_0}$.

  \end{notation}

For Theorem~\ref{f56}, to realize all the tasks, $\lambda > 2^{\aleph_0}$ would suffice;
the requirement in Lemma~\ref{f12.7} that $\lambda=2^\mu$
 is used to get maximal models.  The object of case 3) is to ensure that the
  final model is rich (existentially complete);
    case 4) shows $\uf(M) =\uf(M_*) =\emptyset$.
     After satisfying each task a final section labeled goal verifies that
     each $M_\alpha \in \MM_2$ and so
     $M \in\MM_2$.

    \begin{theorem}\label{f56} There is an $M \in \MM_2$ and in $\bK_1$ that satisfies all the tasks, Thus,
    by Claim~\ref{f44}
  $M \in \bK_2$, and is $P_0$-maximal.
    \end{theorem}


    Proof.  As we construct $M$, we show at appropriate stages that tasks from $\bT_1$ and $\bT_2$
    are satisfied. Further, we show at each stage $\alpha$ the goal: $M_\alpha \in \MM_2$.
 We choose $M_\alpha$ by induction on $\alpha \leq \lambda$ such that:

   \begin{enumerate}
     \item ${\bf w}_\alpha$ witnesses $M_\alpha \in \MM_2$ (Definition~\ref{f50}).
      And for $\beta< \alpha$, $w_\alpha$ extends $w_\beta$.  That is, for $d\in P^{M_\beta}_2$, $\alpha_d[{\bf w}_\alpha] = \alpha_d[{\bf w}_\beta]$, $\eta_d[{\bf w}_\alpha] = \eta_d[{\bf w}_\beta]$, and $a_{d,n}[{\bf w}_\alpha] = a_{d,n}[{\bf w}_\beta]$ . 
     \item $P_2^{M_\alpha} \subseteq P_2^{M_*}$ has cardinality at most $|\alpha| +2^{\aleph_0}$.

    %



         \item if $\alpha = \beta +1$ and $\tbar_\beta$ is relevant to $M_\beta$, $M_\alpha$ satisfies task $\tbar_\beta$.

             \end{enumerate}

\begin{description}
             \item [case 1] If $\alpha = 0$, set $M_0 = M_*\myrestriction (P_0^{M_*} \cup  P_1^{M_*})$.

This condition will be preserved by the induction for all $\alpha$.
              \item [case 2]   Take unions at limits.

             At the successor stage, we now verify task $t_{\beta+1}$ for each of two different types of task.
              Then, we will consider the two cases together
to show the goal  that $M = \bigcup_{\alpha< \lambda} M_{\alpha}  \in \MM_2$.

              \item  [case 3] $\alpha = \beta+1$ and
               say, $\tbar_\beta \in {\bT_1}$, say  $\tbar_\beta = (N_1,N_2)$. (Definition~\ref{f39})

{\bf Choose $M_{\alpha}$:}

                 %
                  %


If $N_1$ is not a subset of  $M_\beta$ then the task is irrelevant and let $M_\alpha =M_\beta$
and $\wbar_\alpha = \wbar_\beta$.
If it is, let $\langle a_\ell\colon \ell< m\rangle$  enumerate $P^{N_2}_2-P^{N_1}_2$
 and
$\langle a'_\ell\colon \ell< m \rangle$ enumerate the first $m$ elements of
$P^{M_*}_2-P^{M_\beta}_2$.
Let $M_\alpha$ extend the $P^{M_\beta}_2$ by
 adding $\langle a'_\ell:\ell<m\rangle$ from $P^{M_*}_2$ to form $P^{M_\alpha}_2$.
 It remains to define the   $\wbar_\alpha$ and  $F^{M_\alpha}_k(a'_\ell)$.

Let $U_\alpha  =\{\delta: (\exists b_\nu \in M_\beta)[\nu(0) = \delta]\}$.
 Clearly $|U_\alpha| \leq |\alpha| + 2^{\aleph_0}$ and
\label{Ualph}\begin{multline*}\ (*)\  \{a_{d,k}\colon k< \omega, d\in P^{M_\beta}_2\} \cup
\{b_\nu\colon (\exists d \in P^{M_\beta}_2)\  \nu\in \Tscr_{\alpha_d}\}
\cup P^{M_*}_{4,1} 
\end{multline*} is included
in the subalgebra of $M_*$ generated by the
$$\{b_\rho \colon \exists \beta \in U_\alpha,
\rho(0) = \beta
\}\cup \{b_{\langle \rangle}\} \cup
P^{M_*}_{4,1}.$$

By induction, since $M_\beta \in \MM_2$ there are witnesses $w_\beta =\langle  \alpha_d,
 \eta_d,a_{d,k}\rangle$
 (formally $\langle  \alpha^\beta_d,\eta^\beta_d, a^\beta_{d,k}\rangle $) for each $d\in
 P^{M_\beta}_2$.
 For the new $a'_\ell$, let $\wbar_\alpha(\ell) =\langle \gamma_\ell,\eta_\ell,0^{M_*}\rangle$
 be chosen with the $\gamma _\ell$ as the first
  $m$ even elements of $\lambda-U_\alpha$ and with $\eta_{\ell}(\wbar_\alpha) =\eta_\ell$ chosen\footnote{In case 3,
  we need choose only a single $\eta_\ell$ for each $\ell< m$. In case 4, we
  choose $2^{\aleph_0}$ distinct $d_\eta$.}
  so that  $\eta_\ell(0) = \gamma_\ell$.
We complete the definition of $M_\alpha$ below by choosing the $F^{M_\alpha}_k$ to
satisfy the task.




%

{\bf Task:} We now verify task $\tbar_{\beta+1}$. 
by showing  in two stages that $N_2$ can be embedded over $N_1$ into $M_\alpha$. First we show there is an
embedding of the Boolean algebras; then we define the $F_k$ on the image to put $M_\alpha$ in $K^1_{\aleph_0}$.
Since $N_2 \in \bK^1 _{<\aleph_0}$, $P^{N_2}_1$ is decomposed as a union of the finite free Boolean
algebras\footnote{While the domain of $N_2\subseteq \lambda$, the  $N_2$-interpretation any relation
symbols in $\tau$ on ordinals not in the domain of $N_1$ has nothing to with the
 interpretations in $M_*$ or $M_\beta$.}  $\langle B^{N_2}_i: i \geq {n_{*}}^{N_2}\rangle$
 where\footnote{Technically, we are
defining $n^{M_\alpha}_*$. But the value is set once and for all at stage $\alpha$ so
we just call it by the final name.}, writing  $n_*$ for $n^{N_2}_*$,
$N_2$ is freely generated  over  $B^{N_2}_{n_{*}}$ mod $P^{N_2}_4$ by
$\{F^{N_2}_k (f): k \geq n^{N_2}_*,
 f \in P^{N_2}_2\}$.
 Similarly, we decompose $P^{N_1}_1$
by $\langle B^{N_1}_i: i \geq {n_{*}}^{N_1}\rangle$.


Since $N_1 \subseteq M_*$ and $N_1 \subseteq N_2$, for each element $e \in P^{N_1}_1$ and any $s$,
$$P^{M_*}_{4,s}(a) \leftrightarrow  P^{N_1}_{4,s}(a) \leftrightarrow  P^{N_2}_{4,s}(a).$$ So no
atom in $N_2- N_1$ is below any element of $N_1$.


Let $\cbar =\langle c_0, \ldots c_{p-1}\rangle$ enumerate the atoms of $N_2$ with
the $c_i$ for $i<r$ enumerating those in  $N_2-N_1$;
they are all in $\BB^{N_2}_{n_*}$.  We set $c'_i = c_i$ if $r\leq i <p$,
choose any $r$ atoms $c'_i$ from $M_*- N_1$  and by Claim~\ref{f34}, we can find 
a $t$  (depending on all of the $c'_i$)
such that for all $i$ if $\nu(0) =  \gamma_\ell$ and $k>t$, $b_{\nu\myrestriction k} \wedge c'_i = 0$.


Each $e \in \BB^{N_2}_* - (P^{N_1}_1 \cup \cbar)$  is a finite
join of $c_i$.
(Note $P^{N_2}_4$ is an alias of $\BB^{N_2}_*$.)
Recall $\{F^{N_2}_k (f): k \geq n^{N_2}_*,
 f \in P^{N_2}_2\}$
 is the pre-image of
a basis of $P^{N_2}_1/ P^{N_2}_4$.
For $f \in P^{N_2}_2$, each $F^{N_2}_k(\hat f) \wedge b^{N_2}_* = e \leq b^{N_2}_*$.
Now define $h_\beta$ mapping $N_2$ into $M_\alpha$ by
\begin{enumerate}
\item  $h_\beta \myrestriction
P^{N_1}_1$ is the identity
\item  $h_\beta(c_i)$ is $c'_i$.
\item  For $e \in \BB^{N_2}_*- P^{N_2}_{4,1}$,  $h_\beta(e) = e' =\bigvee_{c_i \leq e} c'_i$.

\item The $b_{\eta_i} \myrestriction (t +k)$ for $k\geq n_*$ are independent mod $P^{M_*}_4$;
for $a_\ell$ in $P^{N_2}_2 - P^{N_2}_1$  set
 $$h_\beta(F_k^{N_2}(a_\ell)) = b_{\eta_i\myrestriction(t+ k)\add 0}
 \vartriangle b_{\eta_i\myrestriction(t+ k \add 1}) 
  \vee e' = F^{M_\alpha}(a'_\ell)$$

  where $e' = h_\beta(e)$ and $e = F_k^{N_2}(a_\ell) \wedge b^{N_2}_*$.
%
     \item Since the $F_k^{N_2}(a_\ell)$ freely generate $N_2/N_1$ modulo the atoms, $h_\beta$
     extends to an embedding of $N_2$ into $M_\alpha$.


\end{enumerate}

Check using Claim~\ref{backgrbauf}.3 that step 4) is a homomorphism.

We now show $M_\alpha \in \MM_2$. To clarify notation,  by setting\footnote{The $a_{d,n}$ are dummies in this case to provide uniformity
 with  case 4 in proving Lemma~\ref{f53}.}
 $a_{d_\ell,k}=0$ for $i<m$,  we declared:
   $$F_k^{M_\alpha}(d_\ell) = (b_{\eta_i\myrestriction k\add 0} \vartriangle
    b_{\eta_i\myrestriction k \add 1})\vartriangle a_{d_\ell,k}.$$
By Lemma~\ref{countinc}, for some
$n$, for all $k\geq n$,
 $ a \nleq_{P^{M_*}_1} F_k^{M_\alpha}(d_i)$ so condition \ref{f50}.B.2 holds.

  Finally, applying Remark~\ref{backgrbauf}.1  to $f_i = b_{\eta_i\myrestriction k\add 0}  \vartriangle
    b_{\eta_i\myrestriction k \add 1}$ the $\{F_k^{M_\alpha}(d_i): k^Y_1 \leq k <\omega\}$
    are independent  for each $i$ and form a basis for a subalgebra $N'_2$ of $P^{M_*}_1$ over $N_1$.
 Thus, $N'_2 \in \bK^1_{<\aleph_0}$ and we have verified that
    task $\tbar_{\beta+1}$ is satisfied.

\bigskip

 \item  [case 4] $\alpha = \beta+1$ and $\tbar_\beta \in \bT_2$; say, $\tbar_\beta = c$.







We define $M_\alpha$.
 Define $U_\alpha$ as in  Case 3, but  extending $U_\alpha$ to $U'_\alpha$ by
   adding the ordinal named by $c$ if $c \not \in M_\beta$.
   This extension  guarantees that the $F^{M_*}_k(c)$ are in $\BB^0_Y$.
   Now choose
  an even ordinal  $\gamma$ in $\lambda -U_\alpha$ such
   that
      $$\langle \{b_\eta: \eta(0) = \gamma\}\rangle  \cap  \{b_\eta: \eta(0) \in U'_\alpha\}  = \emptyset.$$

Extend  $P^{M_\beta}_2$ by adding a $d_\eta \in P^{M_*}_2$ to  $P^{M_\alpha}_2$ for each $\eta$ with
$\eta(0) = \gamma$.
 To define $F_k^{M_\alpha}(d_\eta)$, for each $\eta \in \lim \Tscr_\gamma$ and
  $k< \omega$,
 choose $i_0< i_1 \leq 2$ that are different from   $\eta(k)$.
                Recalling
                        $c =\tbar_\beta$,  let
$$F_k^{M_\alpha}(d_\eta) = (b_{\eta\myrestriction k\add i_0}
 \vartriangle b_{\eta\myrestriction k \add i_1})\vartriangle (F_k^{M_*}(c)).$$
 Since $M_* \in \bK_{-1}$ for each $a \in P^{M_*}_1$ for all but finitely many
 $n$, $a \wedge F_k^{M_*}(c) =0$.
Thus, for the  $d \in P^{M_\alpha}_2 - P^{M_\beta}_2 $, chosen towards satisfying
 $\tbar_\beta = c$, we have
set $\langle \alpha_d,\eta^\alpha_d,   a_{d,k}\rangle =\langle \gamma , d_\eta,  F^{M_*}_k(c)\rangle$.
 That is, $a_{d,k} = F^{M_*}_k(c)$. Thus, by Lemma~\ref{countinc} for any atom $a$ and all
 but finitely many $n$, $a \wedge F_k^{M_\alpha}(c) =0$ and the countable incompleteness
 requirement is satisfied.

{\bf Task:} We must show $M_\alpha$
 satisfies task $\tbar_\beta$.
Since $\uf(M_*)= \emptyset$, for any non-principal ultrafilter
$D$, there is an $e\in P^{M_*}_2$ such that  
the set $S^{M_*}_e(D) =\{n\colon F^{M_\alpha}_n(e) \in D\}$ is infinite (Definition~\ref{f5}).
 By the
definition of the task $\tbar_\beta=c$,
there is a $D$
where the given $c$   witnesses for $D$ in $\uf(M_*)$.
 We show task $\tbar_\beta$ is satisfied for $D$  by
one of the $d_\eta$, which thus is a witness
to $D \not \in \uf(M_\alpha)$.

 Define $\eta^D \in \lim(\Tscr_\gamma)$ by induction\footnote{This argument
 is patterned on the simple black box in  Lemma 1.5 of \cite{Sh309}, but
  even simpler.}: $\eta^D(0) = \gamma$.
By
Remark~\ref{backgrbauf}.2 one of the three elements
$b_{\langle \gamma, i\rangle} \triangle b_{\langle \gamma, j\rangle}$,
 for $i \neq j$ and $i,j<
 3$, must not be in $D$. Let $\eta^D (1)$ be the other member of $\{0,1,2\}$.
 For $k\geq 1$, suppose $\nu = \eta^D\myrestriction k$ has been defined.  Again, by
Remark~\ref{backgrbauf}.2 one of the three elements
$b_{\nu\add i} \triangle b_{\nu\add j}$, for $i \neq j$ and $i,j<
 3$, must not be in $D$. Let $\eta^D (k)$ be third of the symmetric differences, which
 by Remark~\ref{backgrbauf}.2.c must be in $D$.
  For the  infinitely many $n$ with    $F^{M_\alpha}_n(c) \in D$, we have
 $F^{M_\alpha}_n(d_{\eta^D}) \in D$.

Now we establish the goal for both cases.

 {\bf Goal: $M_\alpha \in \MM_2$:}
 To show $M \in \bK_{-1}$ (and so in $\MM_1$,
 Definition~\ref{f37})  note that  countable
incompleteness
(Definition~\ref{f1}.vii) in each case separately.
For $M_\alpha \in \MM_2$, we
   show $M_\alpha$ satisfies Definition~\ref{f50}.
   The descriptive portions of Conditions A and B.i) of Definition~\ref{f50}  are
clearly satisfied by the construction; Condition B.ii) was shown in the proof of
each case.

For condition C,
choose any finite $Y \subset  P^{M_\alpha}_2$  and partition $Y$ into
 $Y_1 =Y\cap P^{M_\beta}_2$ and $Y_2 = Y-Y_1$. 
%
Set
 $k_1 =k^1_Y$ as the least integer\footnote{Naturally this is only relevant when
  $\alpha_d = \alpha_e$ but than can happen in case 3
 and must happen in case 4.} such that for all $\eta_d \neq \eta_e$ with $d,e\in Y$,
 $\eta_d\myrestriction k^1 \neq \eta_e\myrestriction k^1 $.

For those $d \in Y_1$, we just leave $\wbar_\alpha = \wbar_\beta$. For $d\in Y_2$,
 the two cases differ slightly.

 In case 3, $d \in P^{M_\alpha}_2 - P^{M_\beta}_2 = Y_2$,
  we (implicitly) defined $\wbar_d(\alpha) =
 \langle \alpha_d,\eta_d,  0\rangle$. In case\footnote{Note that in case 3,
  $a_{d,n}$ is constant, while in case 4 it depends on $n$.} 4
 the elements of $Y_2$ are among the $2^{\aleph_0}$ $d_\eta$ with
 $\eta(0)=\gamma$. For them, $\wbar_d(\alpha) =
 \langle \gamma,\eta_d, F^{M_*}_n(c)\rangle$.

For Condition~\ref{f50}.C, we show every element of $W$ is in the $\langle\{b_\nu; \nu(0) \in U_\alpha\}
\rangle$ and so in $\BB^0  _Y$.  By the first line of the proof of case 4, we only need to consider the
$F^M_i(d)$ for $d\in Y$ and $i<|Y|$.

$\{F^M_i(d_k);i<|Y|, k< k^Y_1\}$.
 For $d \in Y_1$ this follows since $d\in P^{M_\beta}_2$
implies $\alpha_d \in U_\alpha$.
In case 3, for $d\in Y_2$ the $a_{d_i,n}$ are all $0$ and the $F^M_n(d)$ for $i<n$ and $n<\omega$
are all Boolean combinations of elements $b_\nu$ with $\nu \unlhd \eta_i\restriction k_1$.


 The difference for case 4
  is in verifying the $a_{d,n}$ are in $U_\alpha$.  
      Now if $a_{d_\eta,n}[{\bf w}_\alpha] = F^{M_*}_n(c)$ is $b_\nu$ then
 $\nu(0) \in U_\alpha$ by the revised definition of $U_\alpha$ in case 4 and so $a_{d,n} \in \BB^0$.

Now, let $M = \bigcup_{\alpha< \lambda}M_\alpha$. Then, $M\in \MM_2$, $|P^M_2|
 = \lambda$. By Lemma~\ref{f53}, $M \in \bK_1$ and each task has been satisfied, so by
 Claim~\ref{f44},  $M \in \bK_2$.

$\qed_{\ref{f56}}$   \end{description}
%



This yields.

                \begin{conclusion} The $M \in \bK_2$ constructed in Theorem~\ref{f56}
                $P_0$-maximal and  all $|P^M_i| = \lambda$.
As in \cite[Corollary 3.3.14]{BaldwinShelahhanfmax}, for every $\lambda$ less
 than the first measurable, since $M\in \bK_2$ implies $|M| \leq 2^{P^{M}_0}$,
  there is a maximal model $M \in \bK_2$ with $2^\lambda \leq |M| <2^{2^\lambda}$.
                \end{conclusion}

%

\begin{question}
\begin{enumerate}
\item
Is there a $\kappa <\mu$, where $\mu$ is the first measurable, such that if a complete sentence has a maximal model in cardinality $\kappa$, it has maximal models in cardinalities cofinal in $\mu$?

    \item Is there a complete sentence that has maximal models cofinally in some $\kappa$ with $\beth_{\omega_1} < \kappa < \mu$ where $\mu$
        is the first measurable, but no larger models are maximal.  Could the first inaccessible be such a $\kappa$?
        \end{enumerate}

        \end{question}


\begin{thebibliography}{KLH16}

\bibitem[BB17]{BBHanf}
John~T. Baldwin and William Boney.
\newblock Hanf numbers and presentation theorems in {A}{E}{C}.
\newblock In Jose Iovino, editor, {\em Beyond First Order Model Theory}, pages
  81--106. Chapman Hall, 2017.

\bibitem[BKS09]{BKS}
John~T. Baldwin, A.~Kolesnikov, and S.~Shelah.
\newblock The amalgamation spectrum.
\newblock {\em Journal of Symbolic Logic}, 74:914--928, 2009.

\bibitem[BKS16]{BKSoul}
John~T. Baldwin, M.~Koerwien, and I.~Souldatos.
\newblock The joint embedding property and maximal models.
\newblock {\em Archive for Mathematical Logic}, 55:545--565, 2016.

\bibitem[BS19]{BSoul}
John~T. Baldwin and Ioannis Souldatos.
\newblock Complete $\mathcal{L}_{\omega_1,\omega}$ with maximal models in
  multiple cardinalities.
\newblock {\em Mathematical Logic Quarterly}, 65(4):444--452, 12 2019.

\bibitem[BS2x]{BaldwinShelahhanfmax}
John~T. Baldwin and S.~Shelah.
\newblock Hanf numbers for extendibility and related phenomena.
\newblock submitted: Shelah number 1092; first posted 2016;
  \url{http://homepages.math.uic.edu/~jbaldwin/pub/ahanfmaxoct3118.pdf}, 202x.

\bibitem[BU17]{BoUn}
W.~Boney and S.~Unger.
\newblock Large cardinal axioms from tameness in {AEC}s.
\newblock {\em Proceedings of the American Mathematical Society},
  145:4517--4532, 2017.


\bibitem[GS05]{GobelShelah}
R.~G{\"o}bel and S.~Shelah.
\newblock How rigid are reduced products.
\newblock {\em Journal of Pure and Applied Algebra}, 202:230--258, 2005.

\bibitem[Gr{\"a}79]{Gratzer}
George Gr{\"a}tzer.
\newblock {\em Universal Algebra}.
\newblock Springer-Verlag, 1979.


\bibitem[Hjo02]{Hjorthchar}
Greg Hjorth.
\newblock Knight's model, its automorphism group, and characterizing the
  uncountable cardinals.
\newblock {\em Journal of Mathematical Logic}, pages 113--144, 2002.

\bibitem[KLH16]{KLH}
Alexei Kolesnikov and Christopher Lambie-Hanson.
\newblock The {H}anf number for amalgamation of coloring classes.
\newblock {\em Journal of Symbolic Logic}, 81:570--583, 2016.

\bibitem[Mag16]{Magcrm}
M.~Magidor.
\newblock Large cardinals and strong logics: {CRM} tutorial lecture 1.
\newblock \url{http://homepages.math.uic.edu/~jbaldwin/pub//MagidorBarc.pdf},
  2016.

\bibitem[Mor65]{Morley65a}
M.~Morley.
\newblock Omitting classes of elements.
\newblock In Addison, Henkin, and Tarski, editors, {\em The Theory of Models},
  pages 265--273. North-Holland, Amsterdam, 1965.

\bibitem[She]{Sh309}
S.~Shelah.
\newblock Black boxes.
\newblock paper 309 archive.0812.0656.

\bibitem[She13]{Sh932}
S.~Shelah.
\newblock Maximal failures of sequence locality in a.e.c.
\newblock preprint on archive: Sh index 932, 2013.

\end{thebibliography}
\end{document}